\documentclass[3p, preprint]{elsarticle}
\usepackage{ecrc}
\volume{00} % '00' means 'unknown volume'
\firstpage{1} % set the starting page if not 1

\makeatletter
\renewcommand\paragraph{\@startsection{paragraph}{4}{\z@}%
           {12\p@ \@plus 6\p@ \@minus 3\p@}%
           {\p@}%
           {\normalfont\normalsize\itshape}}

\makeatother

\usepackage{ecrc_arxiv1}

\usepackage[ruled,english,onelanguage,linesnumbered,vlined]{algorithm2e} % vlined: no 'end' at ends of scopes
\usepackage{amsfonts}
\usepackage{amsmath}
\usepackage{amssymb}
\usepackage{amstext}
\usepackage[title]{appendix}
\usepackage{array}
\usepackage{cancel}
\usepackage{dsfont}
\usepackage{multirow}

\usepackage{epsfig}
\usepackage{subcaption}
\usepackage{calc}
\usepackage{hyperref}
\usepackage{lineno}
\usepackage{multicol}
\usepackage{multirow}
\usepackage{nicefrac}
\usepackage[dvipsnames]{xcolor} % load it before orcidlink (otherwise, orcidlink loads it with no options and we cannot reload it with dvipsname)
\usepackage{orcidlink}
\usepackage{pslatex}
\usepackage{scalerel}
\usepackage{setspace}
\usepackage{txfonts}

\SetCommentSty{mycommfont}
\SetNoFillComment

\newcommand{\h}{\hspace{5mm}}
\newcommand{\hsd}{\hspace{2.5mm}}

\newcommand{\ite}{\item[$\bullet$]}

\newcommand{\dt}{\delta t}

\newcommand{\dsp}{\displaystyle}

\newcommand{\xkinit}{x_k^{\text{init}}}
\newcommand{\tinit}{t^{[\text{init}]}}
\newcommand{\tend}{t^{[\text{end}]}}

\newcommand{\Nmax}{\ensuremath{N_{\text{max}}}}

\newcommand{\nsys}{n_{sys}}

\newcommand{\nstk}[1]{n_{st,#1}}
\newcommand{\nink}[1]{n_{in,#1}}
\newcommand{\noutk}[1]{n_{out,#1}}
\newcommand{\nsttot}{\nstk{tot}}
\newcommand{\nintot}{\nink{tot}}
\newcommand{\nouttot}{\noutk{tot}}

\newcommand{\rnstk}[1]{\mathbb{R}^{\nstk{#1}}}
\newcommand{\rnink}[1]{\mathbb{R}^{\nink{#1}}}
\newcommand{\rnoutk}[1]{\mathbb{R}^{\noutk{#1}}}
\newcommand{\rnsttot}{\mathbb{R}^{\nstk{tot}}}
\newcommand{\rnintot}{\mathbb{R}^{\nink{tot}}}
\newcommand{\rnouttot}{\mathbb{R}^{\noutk{tot}}}
\newcommand{\Insys}{[\![1, \nsys]\!]}
\newcommand{\Instk}[1]{[\![1, \nstk{#1}]\!]}
\newcommand{\Inink}[1]{[\![1, \nink{#1}]\!]}
\newcommand{\Inoutk}[1]{[\![1, \noutk{#1}]\!]}
\newcommand{\Insttot}{[\![1, \nstk{tot}]\!]}
\newcommand{\Inintot}{[\![1, \nink{tot}]\!]}
\newcommand{\Inouttot}{[\![1, \noutk{tot}]\!]}
\newcommand{\Zn}{[\![0, n]\!]}

\newcommand{\bs}{\begin{scriptsize}}
\newcommand{\es}{\end{scriptsize}}

\newcommand{\ddd}[2]{\dsp{\frac{d#1}{d#2}}}
\newcommand{\dspfrac}[2]{\dsp{\frac{#1}{#2}}}

\newcommand{\up}[1]{\textsuperscript{#1}}

\newcommand{\fornits}{F$ _3$ORNITS }

\newcommand{\GV}{{\mathcal{G}_{\!V}}}
\newcommand{\GD}{{\mathcal{G}_{\!D}}}
\newcommand{\PV}{{\mathcal{P}_{\!V}}}
\newcommand{\PD}{{\mathcal{P}_{\!D}}}
\newcommand{\RV}{{\mathcal{R}_{\!V}}}
\newcommand{\RD}{{\mathcal{R}_{\!D}}}
\newcommand{\GbV}{{\bb{\mathcal{G}}_{\!V}}}
\newcommand{\GbD}{{\bb{\mathcal{G}}_{\!D}}}
\newcommand{\PbV}{{\bb{\mathcal{P}}_{\!V}}}
\newcommand{\PbD}{{\bb{\mathcal{P}}_{\!D}}}
\newcommand{\RbV}{{\bb{\mathcal{R}}_{\!V}}}
\newcommand{\RbD}{{\bb{\mathcal{R}}_{\!D}}}
\newcommand{\hyL}{{\hat{y}_{\!L}\phantom{}}}
\newcommand{\hyC}{{\hat{y}_{\!C}\phantom{}}}
\newcommand{\hdyL}{{\hat{\dot{y}}_{\!L}\phantom{}}}
\newcommand{\hdyC}{{\hat{\dot{y}}_{\!C}\phantom{}}}
\newcommand{\tfC}{{\tilde{f}_{\!C}\phantom{}}}
\newcommand{\AVelem}{{\mathcal{A}_{\!V}}_{\text{elem}}}
\newcommand{\ADelem}{{\mathcal{A}_{\!D}}_{\phantom{.}\text{elem}}}
\newcommand{\Belem}{\mathcal{B}_{\text{elem}}}
\newcommand{\AcVelem}{{\check{\mathcal{A}}_{\!V}\phantom{}}_{\text{elem}}}
\newcommand{\AcDelem}{{\check{\mathcal{A}}_{\!D}\phantom{}}_{\phantom{.}\text{elem}}}
\newcommand{\Bcelem}{\check{\mathcal{B}\phantom{}}_{\text{elem}}}
\newcommand{\AcVelemfs}{{\check{\mathcal{A}}_{\!V}\phantom{}}_{\text{elem, first step}}}
\newcommand{\AcDelemfs}{{\check{\mathcal{A}}_{\!D}\phantom{}}_{\phantom{.}\text{elem, first step}}}
\newcommand{\Bcelemfs}{\check{\mathcal{B}\phantom{}}_{\text{elem, first step}}}
\newcommand{\ib}{\bar{\imath}}
\newcommand{\jb}{\bar{\jmath}}
\newcommand{\sigmab}{\bar{\sigma}}
\newcommand{\xb}{\underline{\smash{x}}\vphantom{x}}
\newcommand{\ub}{\underline{\smash{u}}\vphantom{u}}
\newcommand{\yb}{\underline{\smash{y}}\vphantom{y}}
\newcommand{\bb}[1]{\underline{\smash{#1}}\vphantom{#1}}
\def\matDU#1#2#3{
	\ensuremath{
	\left(
	\let\oldarraystretch\arraystretch
	\renewcommand{\arraystretch}{#1}
	\begin{array}{c}
	#2 \\ #3
	\end{array}
	\let\arraystretch\oldarraystretch
	\right)
}}
\def\matDUbar#1#2#3{
	\ensuremath{
	\left(
	\let\oldarraystretch\arraystretch
	\renewcommand{\arraystretch}{#1}
	\begin{array}{c}
	#2 \\ \hline #3
	\end{array}
	\let\arraystretch\oldarraystretch
	\right)
}}
\def\matDD#1#2#3#4#5{
	\ensuremath{
	\left(
	\let\oldarraystretch\arraystretch
	\renewcommand{\arraystretch}{#1}
	\begin{array}{cc}
	#2 & #3 \\ #4 & #5
	\end{array}
	\let\arraystretch\oldarraystretch
	\right)
}}
\def\matDDbar#1#2#3#4#5{
	\ensuremath{
	\left(
	\let\oldarraystretch\arraystretch
	\renewcommand{\arraystretch}{#1}
	\begin{array}{c|c}
	#2 & #3 \\ \hline #4 & #5
	\end{array}
	\let\arraystretch\oldarraystretch
	\right)
}}

\newcommand{\padUp}{\vphantom{\LARGE A}}
\newcommand{\padDown}{\vspace{0.6mm}}
\newcommand{\PadUp}{\vphantom{\huge A}}
\newcommand{\PadDown}{\vspace{1.5mm}}
\newcommand{\Xid}{\mathrm{\Xi}}

\def\wordwrap#1#2{
\let\oldarraystretch\arraystretch
\renewcommand{\arraystretch}{#1}
\begin{tabular}{c}
#2
\end{tabular}
\let\arraystretch\oldarraystretch
}

\def\lwordwrap#1#2{
\let\oldarraystretch\arraystretch
\renewcommand{\arraystretch}{#1}
\begin{tabular}{l}
#2
\end{tabular}
\let\arraystretch\oldarraystretch
}

\def\slwordwrap#1#2#3{
\let\oldarraystretch\arraystretch
\renewcommand{\arraystretch}{#2}
\begin{#1}
\begin{tabular}{l}
#3
\end{tabular}
\end{#1}
\let\arraystretch\oldarraystretch
}

\newcolumntype{P}[1]{>{$}>{\centering\arraybackslash}p{#1}<{$}}

\newcommand{\figuresdir}{.}

\begin{document}

\begin{frontmatter}

\title{MISSILES: an Efficient Resolution of the Co-simulation Coupling Constraint on Nearly Linear Differential Systems through a Global Linear Formulation\tnotemark[1]}

\tnotetext[1]{
Authors Yohan EGUILLON and Bruno LACABANNE are currently employees at Siemens Digital Industries Software. The authors would like to thank Siemens Digital Industries Software for supporting this work, and Institut Camille Jordan and Universit\'e de Lyon for supervising this research. The patent "Advanced cosimulation scheduler for dynamic system simulation" is currently pending to Siemens, and includes the MISSILES co-simulation method.
}

\author[icj,sisw]{\underline{Yohan Eguillon}\orcidlink{0000-0002-9386-4646}}
\author[sisw]{Bruno Lacabanne\orcidlink{0000-0003-1790-3663}}
\author[icj]{Damien Tromeur-Dervout\orcidlink{0000-0002-0118-8100}}
\address[icj]{Institut Camille Jordan, Université de Lyon, UMR5208 CNRS-U.Lyon1, Villeurbanne, France\\
	\{yohan.eguillon,damien.tromeur-dervout\}@univ-lyon1.fr}
\address[sisw]{Siemens Digital Industries Software, Roanne, France\\
	\{yohan.eguillon,bruno.lacabanne\}@siemens.com}
\runauth{Y. Eguillon et al.}

%\linenumbers

\begin{abstract}
In a co-simulation context, interconnected systems of differential equations are solved separately but they regularly communicate data to one another during these resolutions. Iterative co-simulation methods have been developed in order to enhance both stability and accuracy. Such methods imply that the systems must integrate one or more times per co-simulation step (the interval between two consecutive communications) in order to find the best satisfying interface values for exchanged data (according to a given coupling constraint). This requires that every system involved in the modular model is capable of rollback: the ability to re-integrate a time interval that has already been integrated with different input commands. In a paper previously introduced by Eguillon \textit{et al.} in 2022, the COSTARICA process is presented and consists in replacing the non-rollback-capable systems by an estimator on the non-last integrations of the iterative process. The MISSILES algorithm, introduced in this paper, consists in applying the COSTARICA process on every system of a modular model simulated with the IFOSMONDI-JFM iterative co-simulation method (introduced by Eguillon \textit{et al.} in 2021). Indeed, in this case, the iterative part on the estimators of each system can be avoided as the global resolution on a co-simulation step can be written as a single global linear system to solve. Consequently, MISSILES is a non-iterative method that leads to the same solution than the IFOSMONDI-JFM iterative co-simulation method applied to systems using the COSTARICA process to emulate the rollback.
\end{abstract}


\end{frontmatter}

% for footnotes:
\renewcommand*{\thefootnote}{\fnsymbol{footnote}}
\setcounter{footnote}{0}

\section{Introduction}
\label{section:introduction}

Simulation nowadays has a wide range of facets. Among them all, simulation in time of multiphysics models, industrial modular models and composite equations in general is a challenge that is mainly solved by application of a co-simulation method. A co-simulation consists in simulations of interconnected systems exchanging coupling data to one another, each of them embedding its own solver.

Such approach has many benefits. A multiphysical system can be splitted into different systems (also called subsystems) so that each of them represents the equations of a given physical domain and is solved by a tailored solver. For instance, the electrical part can use a dedicated method based on the known shape of the electrical signals, the fluid part can benefit from a solver preserving the conservation laws, etc. Another aspect of it is that the systems can be black-boxed: indeed, given a specified set of interactions, systems can be part of a co-simulation modular model without any need to disclose it in order to simulate it. Among the related industrial applications, we can mention the protection of the intellectual property around the modelling and simulation technologies. Another one is the interoperability: it is possible to connect a set of systems with other systems coming from different platforms and as far as each system provides the requires interactions, the modular model made of these interconnected black-boxed systems can be simulated thanks to a co-simulation method.

A wide range of co-simulation methods (also called co-simulation algorithms) have been proposed in the literature \cite{Gomes2018survey}. Different approaches can be distinguished: from the most generic ones (preserving the black-box aspect of the systems) \cite{Kubler2000} \cite{Sicklinger2014} \cite{Busch2019} \cite{Eguillon2022F3ornits} to the ones using the knowledge of the quantities inside of the systems, benefiting from this information in order to design a more accurate method \cite{Benedikt2013NEPCE} \cite{Stettinger2014ECC} \cite{Stettinger2014IEEEConfControl} \cite{Sadjina2017} \cite{Sadjina2020}. In addition to the fact that the less a method is generic the more it is accurate in general, another trade-off has to be taken into account: the need for advanced capabilities.

A set of basic interactions is always required on a system involed in a co-simulation: the capability to take into account quantities (system's inputs), the capability to move forward in time (simulate on a given time interval) and the capability to retrieve some quantities (system's outputs). No disclosure of the system is required so far. More advanced interaction exist, and some of them are also generic and non-disclosing: for instance, the capability to provide the directional derivatives or the possibility to move forward in time for periods of different sizes. Even though most of the simulation and modelling platforms produce systems with a specified way to know which interactions is available or not on a system, a standard specification exists: the FMI (functional mock-up interface) \cite{FMIStandard}. In this paper, advanced capabilities of the systems will be referred to as in the FMI standard, without loss of generality (other platforms are expected to provide similar capabilities, should it have a different name).

One of these capabilities in particular is critical in co-simulation: the \textit{rollback}. The rollback is the ability of a system to be integrated more than once on a given time interval, also called co-simulation step (or co-simulation time-step). Once a non-rollback-capable system reaches a given time, it can only move forward from this time: its past is frozen and can not be changed anymore. In opposite, a rollback-capable system can be simulated on its last co-simulation step with different inputs than the ones being used in the previous integrations on this co-simulation step. A particular class of co-simulation method requires the rollback capability on every system: the iterative co-simulation methods. These methods use the results of several integrations on a given co-simulation step in order to converge to a more reliable solution regarding the coupling quantities \cite{Kubler2000} \cite{Arnold2001} \cite{Schweizer2014PredCorr} \cite{Sicklinger2014} \cite{Eguillon2019Ifosmondi} \cite{Eguillon2021IfosmondiJFM}.

The rollback being a rare capability in practice, iterative co-simulation methods can usually only be applied on very academic test cases and most of the industrial modular models cannot benefit from the advantages of the iterative co-simulation methods due to the involved non-rollback-capable systems. Moreover, model-based methods are usually challenging to apply on black-boxed systems when the structure of the circuit inside of the system is hidden (including the equations). In order to avoid being restricted to non-iterative and non-disclosing co-simulation methods on industrial systems, the authors proposed an approach to replace the rollback requirement by the use of an estimator based on more common capabilities \cite{Eguillon2022Costarica}. The idea behind this process, called COSTARICA (for: cautiously obtrusive solution to avoid rollback in iterative co-simulation algorithms) is the following: instead of integrating the systems on a given co-simulation step, in case we do not yet know if this integration will be the last one on the step, an estimator of the outputs at the end of the step is used in replacement of the real integration. Once the iterative co-simulation method predicts the converged solution on this co-simulation step (as it iterated on the estimators), a single real integration of the systems is done with the predicted solution as inputs. Reference \cite{Eguillon2022Costarica} details this process and the underlying estimator. The latter requires advanced capabilities which are less rare in practice than the rollback.

The IFOSMONDI-JFM iterative co-simulation algorithm \cite{Eguillon2021IfosmondiJFM}, introduced by the authors as an evolution of the classical IFOSMONDI method \cite{Eguillon2019Ifosmondi} (also introduced by the authors), can benefit from the COSTARICA process in order to run a co-simulation on non-rollback-capable systems. Moreover, in case every system of a modular model uses COSTARICA, the underlying coupling constraint satisfaction problem that the IFOSMONDI-JFM method solves can be transformed into a global linear problem given the expressions of the COSTARICA estimators on every system. The solution of this global linear system directly gives the expressions of the coupling quantities. This paper introduces MISSILES, the co-simulation methods consisting in writting and solving this linear system are using the obtained expressions of the coupling quantities in the real simulation of the systems on every co-simulation step.

The paper is structured as follows: the formalism is introduced as well as the notations with a recall of the outcomes of the COSTARICA process \cite{Eguillon2022Costarica} and IFOSMONDI-JFM method \cite{Eguillon2021IfosmondiJFM}, then the MISSILES method is introduced and the procedure is detailed both regarding the non-iterative co-simulation method and regarding the practical implementation. Finally, the behavior of MISSILES is shown on two test cases before a discussion about the outcomes and the future of the method to conclude.

\section{Formalism and notations}
\label{section:Formalism_and_notations}

This paper focuses on the cases of modular models made of interconnected ordinary differential equations (ODE) systems. The aim is the application to circuits, also sometimes called $0$D systems, corresponding to simulations in time. The method can be adapted to other cases (DAE, PDE), yet the adaptation is not straightforward and will not be discussed in this paper.

\subsection{Problem dimensions and time-domain}
\label{subsection:problem_dimensions_and_time_domain}

Let $\nsys\in\mathbb{N}^{*}$ denote the number of systems involved in a modular model. Please note that the particular case $\nsys=1$ corresponds to a monolithic simulation where a single system is simulated, without connections to any other system, and where no co-simulation method is required. We can thus suppose that $\nsys\in\mathbb{N}\backslash\{0, 1\}$.

With $k\in\Insys$ being the index of a system, we define the following quantities:

\vspace{-0.7em} %layout
\begin{itemize}
\setlength{\itemsep}{-0.4em} %layout
\ite $\nink{k}\in\mathbb{N}$ the number of input variables of system $(S_k)$
\ite $\noutk{k}\in\mathbb{N}$ the number of output variables of system $(S_k)$
\ite $\nstk{k}\in\mathbb{N}$ the number of state variables of system $(S_k)$
\end{itemize}

Equivalently, the inputs, outputs and states can be seen as vectors of dimension $\nink{k}$, $\noutk{k}$ and $\nstk{k}$ respectively.

Finally, $\tinit\in\mathbb{R}$ and $\tend\in]\tinit, +\infty[$ will respectively denote the start time and end time of the simulation.

\subsection{ODE system and discretization}
\label{subsection:ODE_system_and_discretization}

As stated in the beginning of section \ref{section:Formalism_and_notations}, every system is an ODE system. As we are considering interconnected systems, each of them is sensitive to inputs and has outputs. Let the equation of these system be written the following way:

\vspace{-6mm} %layout

\begin{equation}
\label{eq:ODE_Sk}
\begin{array}{ccc}
	\forall k\in\Insys,\
	(S_k):\left\{
		\begin{array}{lcl}
			  \dspfrac{dx_k}{dt} & = & f_k(t, x_k, u_k) \\
			  y_k & = & g_k(t, x_k, u_k) \\
		\end{array}
	\right.
	&
	\text{ where }
	&
	\begin{array}{l}
		t \in [\tinit, \tend] \\
		x_k \in L([\tinit, \tend], \rnstk{k}) \\
		u_k \in L([\tinit, \tend], \rnink{k}) \\
		y_k \in L([\tinit, \tend], \rnoutk{k}) \\
		x_k(\tinit) = \xkinit \in \rnstk{k} \\
	\end{array}
\end{array}
\end{equation}

In \eqref{eq:ODE_Sk}, $u_k$ represent the inputs, $y_k$ the outputs and $x_k$ the states of system $(S_k)$. These time-dependent quantities may be vectorial in case $\nink{k}>1$, $\noutk{k}>1$ or $\nstk{k}>1$. Hence, an extra subscript will denote the index of the element:

\vspace{-2mm} %layout

\begin{equation}
\label{eq:uyx_double_subscript}
\begin{array}{lclllcll}
	u_k & = & (u_{k, i})_{i\in\Inink{k}} & \text{ with } \forall i\in\Inink{k},\ & u_{k, i} & = & (u_k)_i & \in L([\tinit, \tend], \mathbb{R}) \\
	y_k & = & (y_{k, i})_{i\in\Inoutk{k}} & \text{ with } \forall i\in\Inoutk{k},\ & y_{k, i} & = & (y_k)_i & \in L([\tinit, \tend], \mathbb{R}) \\
	x_k & = & (x_{k, i})_{i\in\Instk{k}} & \text{ with } \forall i\in\Instk{k},\ & x_{k, i} & = & (x_k)_i & \in L([\tinit, \tend], \mathbb{R}) \\
\end{array}
\end{equation}

The $f_k$ and $g_k$ applications in \eqref{eq:ODE_Sk} are respectively called the \textit{derivatives} and the \textit{outputs} functions of the ODE. In co-simulation, the systems are black-boxes with a specific set of interaction among which direct calls to $f_k$ or $g_k$ are usually not possible. Hence, the co-simulation is done through a discretization on a time grid on which each mesh is called a \textit{macro-step} (or \textit{co-simulation step}). The nodes are the \textit{communication times}, and they will be denoted by a superscript ${}^{[N]}$ to avoid the confusion with power exponents. $N$ denotes the time index, and the time domain is partitioned as follows in figure \ref{fig:time_discretization}.

\begin{center}
\includegraphics[scale=0.6]{\figuresdir/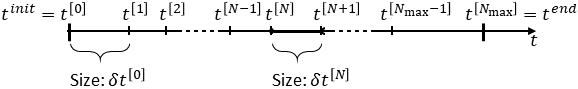}
\vspace{-2mm} %layout
\captionof{figure}{Partition of the time domain in macro-steps}
\label{fig:time_discretization}
\end{center}

Let $\delta t^{[N]}$ denote the size of the $N$\up{th} macro-step:

\vspace{-3mm} %layout

\begin{equation}
\label{eq:macro_step_size}
\forall N \in [\![0, \Nmax-1]\!],\ \delta t^{[N]} = t^{[N+1]}-t^{[N]} > 0
\end{equation}

\vspace{-1mm} %layout

Every system can be integrated on a macro-step as far as an initial condition for the state variables is given. In co-simulation, a system $(S_k)$ can only be integrated on $[t^{[N]}, t^{[N+1]}[$ if it has previously been integrated until $t^{[N]}$. In that case, the initial values for the state variables is $\tilde{x}_k^{[N]}$ the value obtained at this time when integrating the previous step, see \eqref{eq:x_tilde_N}. In the whole paper, quantities with a tilde symbol $\tilde{ }$ will denote an evaluation of a time-dependent quantity at a given time.

\vspace{-7mm} %layout

\begin{equation}
\label{eq:x_tilde_N}
\tilde{x}_k^{[N]} = \left\{
	\renewcommand{\arraystretch}{1.5}
	\begin{array}{ll}
		\xkinit & \text{ if } N=0 \\
		\displaystyle{\lim_{\substack{t \to t^{[N]}\\t<t^{[N]}}}x_k^{[N-1]}(t)} & \text{ otherwise}
	\end{array}
	\renewcommand{\arraystretch}{1.0}
\right.
\end{equation}

\vspace{-1mm} %layout

Equation \eqref{eq:x_tilde_N} define the initial value of the state variables on a time step of the form $[t^{[N]}, t^{[N-1]}[$ with a reference to the time-dependent states $x^{[N-1]}$. The later is the solution of the corresponding system on the corresponding step. Equation \eqref{eq:ODE_Sk_on_tn_tnp} presents these systems.

\vspace{-10mm} %layout

\begin{equation}
\label{eq:ODE_Sk_on_tn_tnp}
\forall k\in\Insys,\
\forall N\in[\![0, \Nmax]\!],\
(S_k^{[N]}):\left\{
	\renewcommand{\arraystretch}{1.5}
	\begin{array}{lcl}
		  \dspfrac{dx_k^{[N]}}{dt} & = & f_k(t, x^{[N]}_k, u^{[N]}_k) \\
		  y_k^{[N]} & = & g_k(t, x_k^{[N]}, u_k^{[N]}) \\
	\end{array}
	\renewcommand{\arraystretch}{1.0}
\right.
\hsd
\text{where}
\hsd
\renewcommand{\arraystretch}{1.3}
\begin{array}{l}
	t \in [t^{[N]}, t^{[N+1]}[ \\
	x_k^{[N]} \in L([t^{[N]}, t^{[N+1]}[, \rnstk{k}) \\
	u_k^{[N]} \in \left(\mathbb{R}_n[t]\right)^{\nink{k}} \\
	y_k^{[N]} \in L([t^{[N]}, t^{[N+1]}[, \rnoutk{k}) \\
	x_k^{[N]}(t^{[N]}) = \tilde{x}^{[N]} \in \rnstk{k} \\
\end{array}
\renewcommand{\arraystretch}{1.0}
\end{equation}

\vspace{-2mm} %layout

Analogously to \eqref{eq:uyx_double_subscript}, the elements of the the inputs, outputs and states on a step use the double subscript notation in this paper:

\vspace{-4mm} %layout

\begin{equation}
\label{eq:uyx_N_double_subscript}
\begin{array}{lclllcll}
	u_k^{[N]} & = & (u_{k, i}^{[N]})_{i\in\Inink{k}} & \text{ with } \forall i\in\Inink{k},\ & u_{k, i}^{[N]} & = & (u_k^{[N]})_i & \in \mathbb{R}_n[t] \\
	y_k^{[N]} & = & (y_{k, i}^{[N]})_{i\in\Inoutk{k}} & \text{ with } \forall i\in\Inoutk{k},\ & y_{k, i}^{[N]} & = & (y_k^{[N]})_i & \in L([t^{[N]}, t^{[N+1]}[, \mathbb{R}) \\
	x_k^{[N]} & = & (x_{k, i}^{[N]})_{i\in\Instk{k}} & \text{ with } \forall i\in\Instk{k},\ & x_{k, i}^{[N]} & = & (x_k^{[N]})_i & \in L([t^{[N]}, t^{[N+1]}[, \mathbb{R}) \\
\end{array}
\end{equation}

\vspace{-1mm} %layout

Please note that, on \eqref{eq:ODE_Sk_on_tn_tnp} and \eqref{eq:uyx_N_double_subscript}, the inputs of $(S_k)$ on the macro-step $[t^{[N]}, t^{[N+1]}[$ are restricted to polynomials. In other words: $u_k^{[N]} \in \left(\mathbb{R}_n[t]\right)^{\nink{k}}$. In most of the co-simulation methods, as the inputs are not known on $[t^{[N]}, t^{[N+1]}[$ when the integration of \eqref{eq:ODE_Sk_on_tn_tnp} is being performed, an extrapolation is made on this interval. Most of the signal extrapolation methods used in co-simulation in practice are covered with the polynomial form: zero-order hold \cite{Sicklinger2014}, first-order hold, Hermite entries \cite{Eguillon2022F3ornits} \cite{Eguillon2021IfosmondiJFM}, smooth polynomial extrapolations \cite{Busch2019}, limited variable order \cite{Eguillon2022F3ornits} \cite{Kraft2019}...

The integration of system $(S_k^{[N]})$ in \eqref{eq:ODE_Sk_on_tn_tnp} leads to the output values $\tilde{y}_k^{[N+1]}$ and their time-derivatives (some co-simulation methods require it) $\tilde{\dot{y}}_k^{[N+1]}$. These quantities are limits at the end of the macro-step, as defined in \eqref{eq:y_tilde_Np1}.

\vspace{-3mm} %layout

\begin{equation}
\label{eq:y_tilde_Np1}
\tilde{y}_k^{[N+1]} = \displaystyle{\lim_{\substack{t \to t^{[N+1]}\\t < t^{[N+1]}}}}y_k^{[N]}(t),
\h\h
\tilde{\dot{y}}_k^{[N+1]} = \displaystyle{\lim_{\substack{t \to t^{[N+1]}\\t < t^{[N+1]}}}}\dspfrac{dy_k^{[N]}}{dt}(t)
\end{equation}

\vspace{-1mm} %layout

As these quantities depend on the initial states on this step $\tilde{x}_k^{[N]}$ and the time-dependent inputs $u_k^{[N]}$, we introduce the formalism of the \textit{step function} $S_k^{[N]}$ (and the analog function $\dot{S}_k^{[N]}$) in \eqref{eq:S_func}. Please note that the step function of $k$\up{th} system on the $N$\up{th} macro-step $S_k^{[N]}$ should not be confused with the system $(S_k^{[N]})$ (with brackets). Nevertheless, the one corresponds to the integration of the other, hence similar notations are not meaningless.

\vspace{-5mm} %layout

\begin{equation}
\label{eq:S_func}
\begin{array}{l}
	S_k^{[N]}:
	\left\{
		\begin{array}{lcl}
			\rnstk{k} \times (\mathbb{R}_n[t])^{\nink{k}} & \rightarrow & \rnoutk{k} \\
			\tilde{x}_k^{[N]}, u_k^{[N]} & \mapsto & S_k^{[N]}(\tilde{x}_k^{[N]}, u_k^{[N]}) = \tilde{y}_k^{[N+1]} \text{ outputs after integration of \eqref{eq:ODE_Sk_on_tn_tnp} at } t^{[N+1]} \text{, see \eqref{eq:y_tilde_Np1}} \\
		\end{array}
	\right.
	\\ \\
	\dot{S}_k^{[N]}:
	\left\{
		\begin{array}{lcl}
			\rnstk{k} \times (\mathbb{R}_n[t])^{\nink{k}} & \rightarrow & \rnoutk{k} \\
			\tilde{x}_k^{[N]}, u_k^{[N]} & \mapsto & \dot{S}_k^{[N]}(\tilde{x}_k^{[N]}, u_k^{[N]}) = \tilde{\dot{y}}_k^{[N+1]} \text{ time-derivatives of these outputs at } t^{[N+1]} \text{, see \eqref{eq:y_tilde_Np1}} \\
		\end{array}
	\right.
\end{array}
\end{equation}

\vspace{-2mm} %layout

An evaluation of the step function $S_k^{[N]}$ implies the integration of the system on the $N$\up{th} macro-step, as shown on figure \ref{fig:Interactions_with_a_system}. The $\dot{S}_k^{[N]}$ function is not available on every system, depending on the modelling and simulation software that generated the system. Similar advanced interactions are not always available and their availability on a given system is often conditioned by capability information. Some of these capabilities are presented in \ref{subsection:Capabilities}.

\subsection{Capabilities}
\label{subsection:Capabilities}

The theoretical framework regarding the systems in a co-simulation has been introduced in subsection \ref{subsection:ODE_system_and_discretization}, yet in practice the black-boxed systems do not always allow to set or evaluate every quantity as we want. Every interaction is standardized in the FMI standard \cite{FMIStandard}, so even if the FMI framework is not mandatory to define, apply or implement the MISSILES method, the standardized naming is useful for the sake of clarity.

Some of the interactions we are interested in in this paper are presented in the table \ref{tab:capabilities}.

\begin{center}
\captionof{table}{Interaction possible with a co-simulation system (non-exhaustive list)}
\label{tab:capabilities}
\renewcommand{\arraystretch}{1.5}
\begin{tabular}{|c|c|c|c|}
	\hline
	Name & Name in the FMI 2.0 standard & Description & Type* \\
	\hline
	\hline
	Provide inputs & \texttt{fmi2SetReal} & \wordwrap{1.2}{Specify $u_k^{[N]}(t^{[N]})$ (zero-order \\hold is used across the \\macro-step by default)} & Basic \\
	\hline
	\wordwrap{1.2}{Provide time-\\dependent inputs} & \texttt{fmi2SetRealInputDerivatives} & \wordwrap{1.2}{Specify non-constant \\$u_k^{[N]}$ (usually polynomial)} & Advanced \\
	\hline
	Do a step & \texttt{fmi2DoStep} & \wordwrap{1.2}{Proceed to the underlying \\integration at the call of \\the $S_k^{[N]}$ function} & Basic \\
	\hline
	Retrieve outputs & \texttt{fmi2GetReal} & Obtain $\tilde{y}_k^{[N+1]}$ & Basic \\
	\hline
	\wordwrap{1.2}{Retrieve outputs \\time-derivatives} & \texttt{fmi2GetRealOutputDerivatives} & Obtain $\tilde{\dot{y}}_k^{[N+1]}$ & Advanced \\
	\hline
	\wordwrap{1.2}{Retrieve state\\values and\\derivatives} & \wordwrap{1.2}{Internal state variables and \\their derivatives must \\be exposed in the system} & \wordwrap{1.2}{\PadUp Obtain $\tilde{x}_k^{[N+1]}$ and $\tilde{f}_k^{[N+1]}$ once \\the macro-step $[t^{[N]}, t^{[N+1]}[$ \\has been integrated} & Advanced \\
	\hline
	Retrieve linearization & \texttt{fmi2GetDirectionalDerivative} & \wordwrap{1.2}{Obtain the derivatives \\of the $f_k$ and $g_k$ functions \\with respect to the states \\and the inputs, respectively} & Advanced \\
	\hline
	Rollback & \wordwrap{1.2}{\texttt{fmi2GetFMUstate}, and \\ \texttt{fmi2SetFMUstate}} & \wordwrap{1.2}{Re-integrate a macro-step \\that has already been integrated} & Advanced \\
	\hline
\end{tabular}
\renewcommand{\arraystretch}{1.0}
\end{center}

*The column named "Type" in table \ref{tab:capabilities} indicates whether the interaction is "Basic", \textit{i.e.} possible in every system for co-simulation, or "Advanced", \textit{i.e.} not mandatory, only available on some systems, depending on the modelling and simulation platform used and submitted to a capability flag (in case of the FMI standard) or a similar mechanism (system generation parameters, for instance).

A visualization of these interactions is presented in figure \ref{fig:Interactions_with_a_system} on a single system (the $k$\up{th} one) and on a macro-step of the form $[t^{[N]}, t^{[N+1]}[$.

The advanced interactions are not equivalently rare in practice: lots of simulation and modelling platforms can generate systems with the possibility to represent polynomial inputs for instance, yet the rollback is very seldom possible on a system for co-simulation. As the rollback is mandatory to use an iterative co-simulation method like \cite{Kubler2000} \cite{Arnold2001} \cite{Bartel2013} \cite{Sicklinger2014} \cite{Eguillon2019Ifosmondi} \cite{Eguillon2021IfosmondiJFM} or any implicit co-simulation method \cite{Schweizer2014PredCorr} \cite{Kraft2021}, lots of industrial models cannot benefit from the advantages of such methods.

However, a slight modification of any co-simulation method requiring the rollback on the systems has been proposed in \cite{Eguillon2022Costarica}: the COSTARICA process.

\begin{center}
\includegraphics[scale=0.27]{\figuresdir/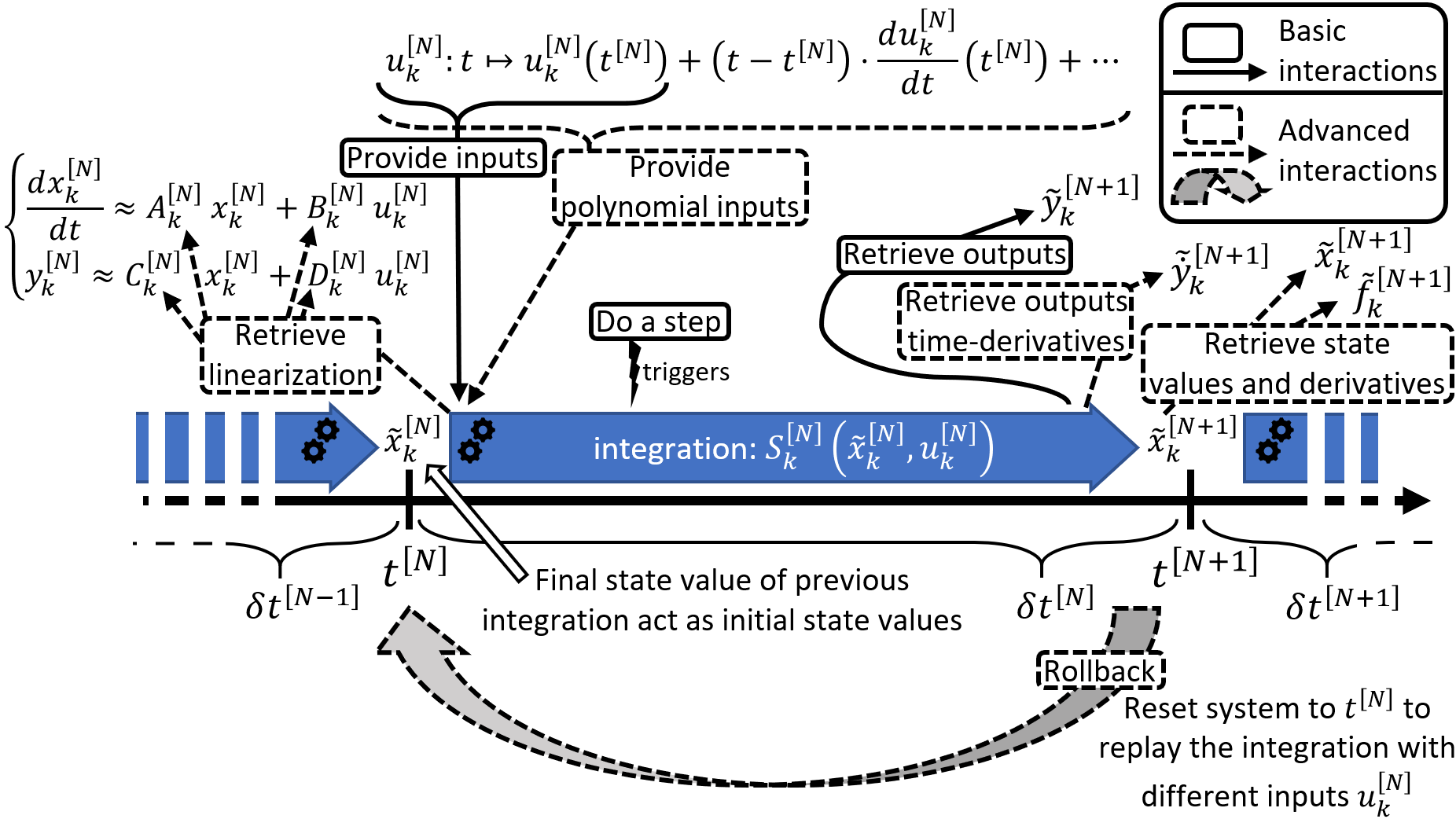}
\vspace{-3mm} %layout
\captionof{figure}{Possible interactions with a system for co-simulation}
\label{fig:Interactions_with_a_system}
\end{center}

\vspace{-3mm} %layout

\subsection{COSTARICA estimator}
\label{subsection:COSTARICA_estimator}

The COSTARICA process, introduced by the authors in \cite{Eguillon2022Costarica}, is a modification applicable to any co-simulation method requiring the rollback as far as:

\vspace{-0.7em} %layout
\begin{itemize}
\setlength{\itemsep}{-0.4em} %layout
\ite the inputs to provide to the systems can be represented by polynomials of a known maximum degree $n\in\mathbb{N}$,
\ite non-rollback-capable systems have the "Retrieve state values and derivatives" and the "Retrieve linearization" capabilities (see table \ref{tab:capabilities}), less rare than the rollback in practice, and
\ite the systems' equations are ODEs (current ongoing developpments tend to extend this restrictions to DAEs, differential algebraic equations).
\end{itemize}

This modification consists in replacing calls to the step function $S_k^{[N]}$ (and $\dot{S}_k^{[N]}$) on every non-rollback-capable system by estimations of it that does not require to integrate the system. By integrating such systems once the definitive inputs for all of them are known, the resulting method does not require the rollback anymore.

The COSTARICA estimators (denoted by a hat symbol $\hat{ }$ ) for output values and time-derivatives are the following ones: regarding a system $(S_k^{[N]})$ (like in \eqref{eq:ODE_Sk_on_tn_tnp}), at a step $[t^{[N]}, t^{[N+1]}[$, the expressions are:

\vspace{-3mm} %layout

\begin{equation}
\label{eq:COSTARICA_y_is_yL_plus_yC}
\renewcommand{\arraystretch}{1.5}
\begin{array}{lclcl}
	\hat{y}_k^{[N+1]} & = & \hyL_k^{[N+1]} + \hyC_k^{[N+1]} & \in & \rnoutk{k} \\
	\hat{\dot{y}}_k^{[N+1]} & = & \hdyL_k^{[N+1]} + \hdyC_k^{[N+1]} & \in & \rnoutk{k} \\
\end{array}
\renewcommand{\arraystretch}{1.0}
\end{equation}

\vspace{-2mm} %layout

\noindent where $\hyC_k^{[N+1]}$ and $\hyC_k^{[N+1]}$ are called the \textit{control parts} and can be obtained by any signal reconstruction method (ZOH, FOH, \fornits \cite{Eguillon2022F3ornits}, ...) as described in \cite{Eguillon2022Costarica}, and where $\hyL_k^{[N+1]}$ and $\hdyL_k^{[N+1]}$ are called the \textit{linear parts} and can be computed with:

\vspace{-7mm} %layout

\begin{equation}
\label{eq:COSTARICA_yL_is_GXi_plus_Pxi}
\renewcommand{\arraystretch}{1.5}
\begin{array}{lcl}
	\hyL_k^{[N+1]} & = & \GV_k^{[N]} \check{\Xid}_k^{[N]} + \PV_k^{[N]} \tilde{x}_k^{[N]} + \RV_k^{[N]} \tfC_k^{[N]} \\
	\hdyL_k^{[N+1]} & = & \GD_k^{[N]} \check{\Xid}_k^{[N]} + \PD_k^{[N]} \tilde{x}_k^{[N]} + \RD_k^{[N]} \tfC_k^{[N]} \\
\end{array}
\renewcommand{\arraystretch}{1.0}
\end{equation}

\vspace{-2mm} %layout

\noindent where $\GV_k^{[N]}$ and $\GD_k^{[N]}$ are tensors of size $\noutk{k}\times\nink{k}\times (n+1)$, where $\PV_k^{[N]}$, $\PD_k^{[N]}$, $\RV_k^{[N]}$ and $\RD_k^{[N]}$ are matrices of size $\noutk{k}\times\nstk{k}$, and where $\tfC_k^{[N]}$ is a vector of size $\nstk{k}$. All seven depend on the current macro-step size $\dt^{[N]}$, on the values at $t^{[N]}$ of the inputs, states and state derivatives, and the linearization of the system $(S_k^{[N]})$ at time $t^{[N]}$ in the form of the matrices referred to as $A^{[N]}$, $B^{[N]}$, $C^{[N]}$ and $D^{[N]}$ on figure \ref{fig:Interactions_with_a_system}. Their expressions are given in \eqref{eq:COSTARICA_GV_GD_PV_PD_RV_RD} \eqref{eq:COSTARICA_gamma_pi_theta} \eqref{eq:COSTARICA_Gamma_Pi_Theta}, and \eqref{eq:COSTARICA_tfC} yet further computation details are given in \cite{Eguillon2022Costarica}.

\vspace{-2mm} %layout

\begin{equation}
\label{eq:COSTARICA_GV_GD_PV_PD_RV_RD}
\begin{array}{ccccc}
	\GV_k^{[N]}	= \gamma_k^{[N]}(\dt^{[N]}),
	& \h &
	\PV_k^{[N]}	= \varpi_k^{[N]}(\dt^{[N]}),
	& \h &
	\RV_k^{[N]}	= \vartheta_k^{[N]}(\dt^{[N]}),
	\\
	\GD_k^{[N]}	= \dfrac{d\gamma_k^{[N]}}{dt}(\dt^{[N]}),
	& \h &
	\PD_k^{[N]}	= \dfrac{d\varpi_k^{[N]}}{dt}(\dt^{[N]})
	& \h &
	\RD_k^{[N]}	= \dfrac{d\vartheta_k^{[N]}}{dt}(\dt^{[N]})
\end{array}
\end{equation}

For all system $k$, at every time step $[t^{[N]}, t^{[N+1]}[$, the $\gamma_k^{[N]}$, $\varpi_k^{[N]}$ and $\vartheta_k^{[N]}$ functions are defined by:

\vspace{-2mm} %layout

\begin{equation}
\label{eq:COSTARICA_gamma_pi_theta}
\gamma_k^{[N]}:\check{t}
\mapsto
\mathcal{L}^{-1}
\left(
	\Gamma_k^{[N]}
\right)(\ \check{t}\ ),
\h\h
\varpi_k^{[N]}:\check{t}
\mapsto
\mathcal{L}^{-1}
\left(
	\Pi_k^{[N]}
\right)(\ \check{t}\ ),
\h\h
\vartheta_k^{[N]}:\check{t}
\mapsto
\mathcal{L}^{-1}
\left(
	\Theta_k^{[N]}
\right)(\ \check{t}\ )
\end{equation}

\vspace{-2mm} %layout

\noindent with

\vspace{-9mm} %layout

\begin{equation}
\label{eq:COSTARICA_Gamma_Pi_Theta}
\Gamma_k^{[N]}:s
\mapsto
\left(
	C^{[N]} (sI-A^{[N]})^{-1}B^{[N]} + D^{[N]}
\right)
\otimes
\left(
\left(
	\frac{p!}{s^{p+1}}
\right)_{p\in\Zn}
\right)^T,
\hsd
\Pi_k^{[N]}:s
\mapsto
C^{[N]} (sI-A^{[N]})^{-1},
\hsd
\Theta_k^{[N]}:s
\mapsto
\frac{1}{s} \Pi_k^{[N]}(s)
\end{equation}

\vspace{-3mm} %layout

Regarding the column vector $\tfC_k^{[N]}$, it is given by:

\vspace{-5mm} %layout

\begin{equation}
\label{eq:COSTARICA_tfC}
\tfC_k^{[N]}
=
\lim_{\substack{t\to t^{[N]}\\t<t^{[N]}}}
\left(
	f_k\big(t, x_k^{[N-1]}(t), u_k^{[N-1]}(t)\big) - A^{[N]} x_k^{[N-1]}(t) - B^{[N]} u_k^{[N-1]}(t)
\right)
=
\tilde{f}_k^{[N]}
- A^{[N]} \tilde{x}_k^{[N]}
- B^{[N]} \lim_{\substack{t\to t^{[N]}\\t<t^{[N]}}}
\left(
	u_k^{[N-1]}(t)
\right)
\end{equation}

\vspace{-3mm} %layout

In \eqref{eq:COSTARICA_yL_is_GXi_plus_Pxi}, $\check{\Xid}_k^{[N]}$ is a matrix of size $\nink{k}\times (n+1)$ and corresponds to the coefficients of the time-shifted polynomials of the input variables $u_k^{[N]}$. On the whole paper, the polynomial degree will be the only dimension in the matrices and tensors whose indexing starts by zero (instead of $1$, like the system index for instance), so that the $p$\up{th} element correspond to the monomial of degree $p$ with $p\in\Zn$. The coefficients of the polynomials of the input variables $u_k^{[N]}$ can be written as:

\vspace{-6mm} %layout

\begin{equation}
\label{eq:Xi_and_ukN_coefficients}
\Xid_k^{[N]} = \left(a_{j, p}^{[N]}\right)_{\substack{j\in\Inink{k}\\p\in\Zn}}
\ \text{where} \
\forall j \in\Inink{k},\ u_{k, j}^{[N]} = \left(u_k^{[N]}\right)_j = \left(t\mapsto\sum_{p=0}^{n}a_{j, p}^{[N]}t^p\right) \in \mathbb{R}_n[t]
\end{equation}

\vspace{-2mm} %layout

\noindent as the inputs are considered polynomial, with a maximum degree of $n$. However, the $\check{\Xid}_k^{[N]}$ required in \eqref{eq:COSTARICA_yL_is_GXi_plus_Pxi} for the COSTARICA estimators is slightly different than $\Xid$, due to a time-shift.

\subsection{Time shift}
\label{subsection:time_shift}

Let's denote by $\check{u}_k^{[N]}$ the time-shifted inputs of $(S_k^{[N]})$ on $[0, \dt^{[N]}[$ (instead of $[t^{[N]}, t^{[N+1]}[$). The time-shifted time variable in $[0, \dt^{[N]}[$ is denoted by $\check{t}=t-t^{[N]}$.

\vspace{-4mm} %layout

\begin{equation}
\label{eq:u_shifted}
\forall \check{t}\in[0, \dt^{[N]}[,\ \check{u}_k^{[N]}(\check{t}) = u_k^{[N]}(t^{[N]} + \check{t})
\end{equation}

\vspace{-2mm} %layout

Let's denote the coefficients of $\check{u}_k^{[N]}$ as follows:

\vspace{-5mm} %layout

\begin{equation}
\label{eq:ukN_check_coefficients}
\forall j \in\Inink{k},\ \check{u}_{k, j}^{[N]} = \left(\check{u}_k^{[N]}\right)_j = \left(\check{t}\mapsto\sum_{p=0}^{n}\check{a}_{j, p}^{[N]}\check{t}^p\right) \in \mathbb{R}_n[t]
\end{equation}

\vspace{-2mm} %layout

Finally, the $\check{\Xid}_k^{[N]}$ matrix in \eqref{eq:COSTARICA_yL_is_GXi_plus_Pxi} is filled with the coefficients of the time-shifted polynomials of the input variables.

\vspace{-3mm} %layout

\begin{equation}
\label{eq:Xi_check}
\check{\Xid}_k^{[N]} = \left(\check{a}_{j, p}^{[N]}\right)_{\substack{j\in\Inink{k}\\p\in\Zn}}
\end{equation}

\vspace{-1mm} %layout

These coefficients can be computed using formula \eqref{eq:a_check_def} proven in \cite{Eguillon2022Costarica}.

\vspace{-3mm} %layout

\begin{equation}
\label{eq:a_check_def}
\forall j\in\Inink{k},\ \forall p\in\Zn,\ \check{a}_{j, p}^{[N]} = \sum_{q=p}^{n} a_{j, q}^{[N]} \left(\begin{smallmatrix}q\\p\end{smallmatrix}\right) (t^{[N]})^{q-p}
\end{equation}

\vspace{-3mm} %layout

\noindent where $\left(\begin{smallmatrix}q\\p\end{smallmatrix}\right)$ denote the \textit{binomial coefficient} of $q$ choose $p$, also referred to as the \textit{combinations} in the literature.

This transformation can be written as a tensor of order $4$ of size $\nink{k}\times (n+1)\times\nink{k}\times (n+1)$, denoted by $\mathcal{C}_k^{[N]}$ and satisfying:

\vspace{-6mm} %layout

\begin{equation}
\label{eq:Xi_check_is_CkN_times_Xi}
\check{\Xid}_k^{[N]} = \mathcal{C}_k^{[N]} \Xid_k^{[N]}
\end{equation}

\vspace{-3mm} %layout

Such a tensor is defined by:

\vspace{-6mm} %layout

\begin{equation}
\label{eq:C_kN}
\mathcal{C}_k^{[N]} = \left(
	\delta_{j_1, j_2}
	\cdot
	\mathds{1}_{q\geqslant p}
	\cdot
	\left(\begin{smallmatrix}q\\p\end{smallmatrix}\right) (t^{[N]})^{q-p}
\right)_{\substack{
	j_1\in\Inink{k} \\
	p\in\Zn \\
	j_2\in\Inink{k} \\
	q\in\Zn \\
}}
\end{equation}

\vspace{-1mm} %layout

\noindent where $\delta_{j_1, j_2}$ denotes the Kronecker coefficient: $1$ in case $j_1=j_2$, and $0$ otherwise. As the coefficients of a single time-shifted input of index $j$ only depends on the original input of index $j$ and none other inputs, $\mathcal{C}_k^{[N]}$ is a \textit{diagonal} tensor with respect to its first and third dimensions. This is precisely the role of the Kronecker coefficient in \eqref{eq:C_kN}, and the diagonal elements with respect to these dimensions is the following matrix:

\vspace{-3mm} %layout

\begin{equation}
\label{eq:C_kN_diag}
\forall j\in\Inink{k},\
\left((\mathcal{C}_k^{[N]})_{j, p, j, q}\right)_{\substack{
	p\in\Zn \\
	q\in\Zn \\
}}
=
\left(
\begin{array}{cccc}
	\left(\begin{smallmatrix}0\\0\end{smallmatrix}\right)(t^{[N]})^{0-0}
	&
	\left(\begin{smallmatrix}1\\0\end{smallmatrix}\right)(t^{[N]})^{1-0}
	&
	\cdots
	&
	\left(\begin{smallmatrix}n\\0\end{smallmatrix}\right)(t^{[N]})^{n-0}
	\\
	0
	&
	\left(\begin{smallmatrix}1\\1\end{smallmatrix}\right)(t^{[N]})^{1-1}
	&
	\cdots
	&
	\left(\begin{smallmatrix}n\\1\end{smallmatrix}\right)(t^{[N]})^{n-1}
	\\
	\vdots & \ddots & \ddots & \vdots \\
	0 & 0 & \cdots & \left(\begin{smallmatrix}n\\n\end{smallmatrix}\right)(t^{[N]})^{n-n}
	\\
\end{array}
\right)
\end{equation}

\subsection{Hermite interpolation}
\label{subsection:hermite_interpolation}

At some point in the method, a Hermite interpolation will be required by the IFOSMONDI-JFM method \cite{Eguillon2021IfosmondiJFM}). The interpolation has to occur on two points in time, with values and first-order time-derivative constraints. The aim of this subsection is to give a linear relationship between the value and time-derivative contraints on one of the points and the coefficients of the Hermite interpolating polynomial.

Let's consider the one-dimensional Hermite interpolation on $2$ generic points at times $t_1$ and $t_2$ (with $t_2 \neq t_1$), with values $v_1$ and $v_2$ and derivatives $\dot{v}_1$ and $\dot{v}_2$ respectively:

\vspace{-5mm} %layout

\begin{equation}
\label{eq:Hermite_2points_generic_properties}
\renewcommand{\arraystretch}{2.0}
\begin{array}{c|c|c}
	\multirow{2}{*}
	{$
		\pi:\left\{
		\renewcommand{\arraystretch}{1.0}
		\begin{array}{lcl}
			\mathbb{R} & \rightarrow & \mathbb{R} \\
			t & \mapsto	& \pi(t) \\
		\end{array}
		\right.
	$}
	& \pi(t_1) = v_1 & \ddd{\pi}{t}(t_1) = \dot{v}_1 \\
	& \pi(t_2) = v_2 & \ddd{\pi}{t}(t_2) = \dot{v}_2 \\
\end{array}
\renewcommand{\arraystretch}{1.0}
\end{equation}

\vspace{-3mm} %layout

Such polynomial has the following expression:

\vspace{-3mm} %layout

\begin{equation}
\label{eq:Hermite_2points_generic_expression}
\pi(t)
=
a_3 t^3 + a_2 t^2 + a_1 t + a_0
\end{equation}

\vspace{-1mm} %layout

\noindent where the coefficients of $\pi$ have the following linear expressions with respect to $v_2$ and $\dot{v}_2$ constraints:

\vspace{-3mm} %layout

\begin{equation}
\label{eq:Hermite_2points_a3_a2_a1_a0}
\begin{array}{lcl}
	a_3
	& = &
	\frac{1}{(t_2-t_1)^2}
	\left(
	\!\!
	\begin{array}{cc}
		\frac{-2}{t_2 - t_1}
		&
		1
	\end{array}
	\!\!
	\right)
	\left(
	\!\!
	\begin{array}{c}
		v_2 \\
		\dot{v}_2
	\end{array}
	\!\!
	\right)
	+
	\frac{1}{(t_2-t_1)^2}
	\left(
	\!\!
	\begin{array}{c}
		\dot{v}_1 + \frac{2v_1}{t_2-t_1}
	\end{array}
	\!\!
	\right)
	\\
	a_2
	& = &
	\frac{1}{(t_2-t_1)^2}
	\left(
	\!\!
	\begin{array}{cc}
		1 + \frac{2(t_2+2t_1)}{t_2-t_1}
		&
		-(t_2+2t_1)
	\end{array}
	\!\!
	\right)
	\left(
	\!\!
	\begin{array}{c}
		v_2 \\
		\dot{v}_2
	\end{array}
	\!\!
	\right)
	+
	\frac{1}{(t_2-t_1)^2}
	\left(
	\!\!
	\begin{array}{c}
		v_1-(t_1+2t_2)(\dot{v}_1+\frac{2v_1}{t_2-t_1})
	\end{array}
	\!\!
	\right)
	\\
	a_1
	& = &
	\frac{1}{(t_2-t_1)^2}
	\left(
	\!\!
	\begin{array}{c}
		t_1\left(\frac{-4t_2}{t_2-t_1}-2-\frac{2t_1}{t_2-t_1}\right)
		\\
		t_1\left(2+t_1\right)
	\end{array}
	\!\!
	\right)^T
	\left(
	\!\!
	\begin{array}{c}
		v_2 \\
		\dot{v}_2
	\end{array}
	\!\!
	\right)
	+
	\frac{1}{(t_2-t_1)^2}
	\left(
	\!\!
	\begin{array}{c}
		t_2\left(2\big((\dot{v}_1+\frac{2v_1}{t_2-t_1})t_1-v_1\big)+t_2(\dot{v}_1+\frac{2v_1}{t_2-t_1})\right)
	\end{array}
	\!\!
	\right)
	\\
	a_0
	& = &
	\frac{1}{(t_2-t_1)^2}
	\left(
	\!\!
	\begin{array}{cc}
		t_1^2
		(
			1
			+
			\frac{2t_2}{t_2-t_1}
		)
		&
		-t_1^2 t_2
	\end{array}
	\!\!
	\right)
	\left(
	\!\!
	\begin{array}{c}
		v_2 \\
		\dot{v}_2
	\end{array}
	\!\!
	\right)
	+
	\frac{1}{(t_2-t_1)^2}
	\left(
	\!\!
	\begin{array}{c}
		t_2^2\left(v_1-(\dot{v}_1+\frac{2v_1}{t_2-t_1})t_1\right)
	\end{array}
	\!\!
	\right)
\end{array}
\end{equation}

\vspace{-1mm} %layout

Finally, we can write:

\vspace{-6mm} %layout

\begin{equation}
\label{eq:Hermite_linexp_general}
\left(
	\begin{array}{c}
		a_0 \\
		a_1 \\
		a_2 \\
		a_3 \\
	\end{array}
\right)
=
\left(
	\AVelem | \ADelem
\right)
\left(
\!\!
\begin{array}{c}
	v_2 \\
	\dot{v}_2
\end{array}
\!\!
\right)
+
\Belem
\end{equation}

\vspace{-1mm} %layout

\noindent where $\left(\AVelem | \ADelem\right)$ is the $4\times 2$ matrix concatenation of the $4\times 1$ column vectors $\AVelem$ and $\ADelem$, and where $\AVelem$, $\ADelem$ and the third $4\times 1$ column vector $\Belem$ have the following expressions:

\vspace{-1mm} %layout

\begin{equation}
\label{eq:Hermite_AVelem_ADelem_Belem}
\begin{array}{c}
	\AVelem
		=
		\frac{1}{(t_2\!-\!t_1)^2}
		\!
		\left(
		\!\!\!
		\renewcommand{\arraystretch}{1.5}
		\begin{array}{c}
			t_1^2
			(
				1
				+
				\frac{2t_2}{t_2-t_1}
			)
			\\
			t_1\left(\frac{-4t_2}{t_2-t_1}-2-\frac{2t_1}{t_2-t_1}\right)
			\\
			1 + \frac{2(t_2+2t_1)}{t_2-t_1}
			\\
			\frac{-2}{t_2 - t_1}
		\end{array}
		\renewcommand{\arraystretch}{1.0}
		\!\!\!
		\right)
	,\h\h
	\ADelem
		=
		\frac{1}{(t_2\!-\!t_1)^2}
		\!
		\left(
		\!\!\!
		\begin{array}{c}
			-t_1^2 t_2
			\\
			t_1\left(2+t_1\right)
			\\
			-(t_2+2t_1)
			\\
			1
		\end{array}
		\!\!\!
		\right)
	\\
	\Belem
	=
	\frac{1}{(t_2\!-\!t_1)^2}
	\left(
	\!\!\!
	\begin{array}{c}
		t_2^2\left(v_1-(\dot{v}_1+\frac{2v_1}{t_2-t_1})t_1\right)
		\\
		t_2\left(2\big((\dot{v}_1+\frac{2v_1}{t_2-t_1})t_1-v_1\big)+t_2(\dot{v}_1+\frac{2v_1}{t_2-t_1})\right)
		\\
		v_1-(t_1+2t_2)(\dot{v}_1+\frac{2v_1}{t_2-t_1})
		\\
		\dot{v}_1 + \frac{2v_1}{t_2-t_1}
	\end{array}
	\!\!\!
	\right)
\end{array}
\end{equation}

\subsection{Connecting systems into a modular model}
\label{subsection:connecting_systems_into_a_modular_model}

Let $\nintot$, $\nouttot$ and $\nsttot$ respectively denote the total amount os inputs, outputs and states:

\vspace{-3mm} %layout

\begin{equation}
\label{eq:sizes_tot}
\nintot = \sum_{k=1}^{\nsys}\nink{k},
\h
\nouttot = \sum_{k=1}^{\nsys}\noutk{k},
\h
\nsttot = \sum_{k=1}^{\nsys}\nstk{k}
\end{equation}

\vspace{-1mm} %layout

Let also the total input, total output and total state vectors be the concatenation of all inputs, outputs and states (respectively) of all systems on a step of the form $[t^{[N]}, t^{[N+1]}[$:

\begin{equation}
\label{eq:ios_tot}
\renewcommand{\arraystretch}{1.5}
\begin{array}{c}
	\ub^{[N]} = (u_k^{[N]})_{k\in\Insys} \in \mathbb{R}^{\nintot},
	\h
	\yb^{[N]} = (y_k^{[N]})_{k\in\Insys} \in L([t^{[N]}, t^{[N+1]}[, \rnouttot),
	\\
	\xb^{[N]} = (x_k^{[N]})_{k\in\Insys} \in L([t^{[N]}, t^{[N+1]}[, \rnsttot)
\end{array}
\renewcommand{\arraystretch}{1.0}
\end{equation}

At any time $t$ in any step $[t^{[N]}, t^{[N+1]}[$, the total input, total output and total state vectors are big column vectors:

\vspace{-6mm} %layout

\begin{equation}
\label{eq:ios_tot_fcts}
\renewcommand{\arraystretch}{1.7}
\begin{array}{lccccr}
	\ub^{[N]}:t\mapsto
	\Big(
		&
		u_{1, 1}^{[N]}(t)\ ,\ ...\ ,\ u_{1, \nink{1}}^{[N]}(t)\ ,
		\hsd & \hsd
		u_{2, 1}^{[N]}(t)\ ,\ ...\ ,\ u_{2, \nink{2}}^{[N]}(t)\ ,
		\hsd & \hsd
		...\ ,
		\hsd & \hsd
		u_{\nsys, 1}^{[N]}(t)\ ,\ ...\ , u_{\nsys, \nink{\nsys}}^{[N]}(t)
		&
	\Big)^T
	\\
	\yb^{[N]}:t\mapsto
	\Big(
		&
		y_{1, 1}^{[N]}(t)\ ,\ ...\ ,\ y_{1, \noutk{1}}^{[N]}(t)\ ,
		\hsd & \hsd
		y_{2, 1}^{[N]}(t)\ ,\ ...\ ,\ y_{2, \noutk{2}}^{[N]}(t)\ ,
		\hsd & \hsd
		...\ ,
		\hsd & \hsd
		y_{\nsys, 1}^{[N]}(t)\ ,\ ...\ , y_{\nsys, \noutk{\nsys}}^{[N]}(t)
		&
	\Big)^T
	\\
	\xb^{[N]}:t\mapsto
	\Big(
		&
		x_{1, 1}^{[N]}(t)\ ,\ ...\ ,\ y_{1, \nstk{1}}^{[N]}(t)\ ,
		\hsd & \hsd
		x_{2, 1}^{[N]}(t)\ ,\ ...\ ,\ y_{2, \nstk{2}}^{[N]}(t)\ ,
		\hsd & \hsd
		...\ ,
		\hsd & \hsd
		x_{\nsys, 1}^{[N]}(t)\ ,\ ...\ , y_{\nsys, \nstk{\nsys}}^{[N]}(t)
		&
	\Big)^T
	\\
\end{array}
\renewcommand{\arraystretch}{1.0}
\end{equation}

\vspace{-2mm} %layout

Let $\tilde{\xb}^{[N]}$, $\tilde{\yb}^{[N]}$ and $\tilde{\dot{\yb}}^{[N]}$ be the instant total state, instant total output and instant total output time-derivative vectors (respectively). Based on \eqref{eq:x_tilde_N} and \eqref{eq:y_tilde_Np1}, we can write, for any $N\in[\![0, \Nmax]\!]$:

\vspace{-2mm} %layout

\begin{equation}
\label{eq:os_tot_snapshots}
\renewcommand{\arraystretch}{1.7}
\begin{array}{lccccr}
	\tilde{\xb}^{[N]} =
	\Big(
		&
		\tilde{x}_{1, 1}^{[N]}\ ,\ ...\ ,\ \tilde{x}_{1, \nstk{1}}^{[N]}\ ,
		\hsd & \hsd
		\tilde{x}_{2, 1}^{[N]}\ ,\ ...\ ,\ \tilde{x}_{2, \nstk{2}}^{[N]}\ ,
		\hsd & \hsd
		...\ ,
		\hsd & \hsd
		\tilde{x}_{\nsys, 1}^{[N]}\ ,\ ...\ ,\ \tilde{x}_{\nsys, \nstk{\nsys}}^{[N]}
		&
	\Big)^T
	\\
	\tilde{\yb}^{[N]} =
	\Big(
		&
		\tilde{y}_{1, 1}^{[N]}\ ,\ ...\ ,\ \tilde{y}_{1, \noutk{1}}^{[N]}\ ,
		\hsd & \hsd
		\tilde{y}_{2, 1}^{[N]}\ ,\ ...\ ,\ \tilde{y}_{2, \noutk{2}}^{[N]}\ ,
		\hsd & \hsd
		...\ ,
		\hsd & \hsd
		\tilde{y}_{\nsys, 1}^{[N]}\ ,\ ...\ ,\ \tilde{y}_{\nsys, \noutk{\nsys}}^{[N]}
		&
	\Big)^T
	\\
	\tilde{\dot{\yb}}^{[N]} =
	\Big(
		&
		\tilde{\dot{y}}_{1, 1}^{[N]}\ ,\ ...\ ,\ \tilde{\dot{y}}_{1, \noutk{1}}^{[N]}\ ,
		\hsd & \hsd
		\tilde{\dot{y}}_{2, 1}^{[N]}\ ,\ ...\ ,\ \tilde{\dot{y}}_{2, \noutk{2}}^{[N]}\ ,
		\hsd & \hsd
		...\ ,
		\hsd & \hsd
		\tilde{\dot{y}}_{\nsys, 1}^{[N]}\ ,\ ...\ ,\ \tilde{\dot{y}}_{\nsys, \noutk{\nsys}}^{[N]}
		&
	\Big)^T
	\\
\end{array}
\renewcommand{\arraystretch}{1.0}
\end{equation}

\vspace{-1mm} %layout

The modular model is the global co-simulation model made of the interconnected $\nsys$ systems. The outputs of the different systems are connected to inputs of other systems, with a possibility to connect a single output to several inputs. An simple example with $\nsys=3$ systems is presented on figure \ref{fig:connections_phi}.

The connections between the systems will be denoted by a matrix filled with zeros and ones, with $\nouttot$ rows and $\nintot$ columns denoted by $\Phi$.

\vspace{-4mm} %layout

\begin{equation}
\label{eq:Phi}
\forall \ib\in\nouttot,\ \forall \jb\in\nintot,\ 
\Phi_{\ib, \jb} = \left\{
\begin{array}{ll}
	1 & \text{if output $\ib$ is connected to input $\jb$} \\
	0 & \text{otherwise}
\end{array}
\right.
\end{equation}

\vspace{-1mm} %layout

Please note that if each output is connected to exactely one input, $\Phi$ is a square matrix. Moreover, it is a permutation matrix in this case. Otherwise, if an output is connected to several inputs, more than one $1$ appears at the corresponding row of $\Phi$. Without loss of generality, let's consider that there can neither be more nor less than one $1$ on each column of $\Phi$ considering that an input can neither be connected to none nor several outputs. Indeed, a system with an input connected to nothing is not possible (a value has to be given), and a connection of several outputs in the same input can always be decomposed regarding a relation (sum, difference, ...) so that this situation is similar to distinct inputs connected to a single output each, with these inputs being combined (added, substracted, ...) inside of the considered system.

\vspace{-3mm} %layout

\begin{center}
\includegraphics[scale=0.39]{\figuresdir/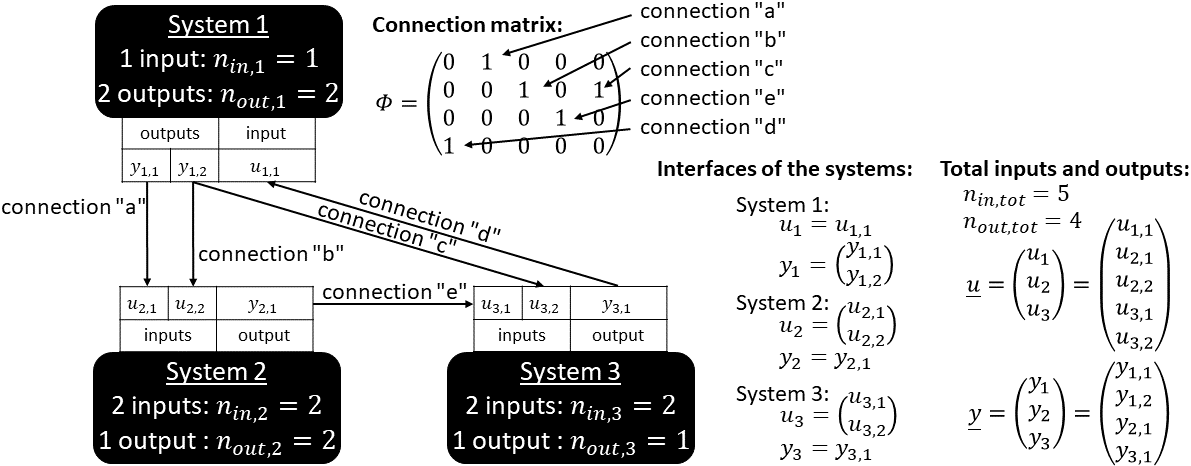}
\vspace{-3mm} %layout
\captionof{figure}{Example of a $3$-system co-simulation model with its interfaces and its $\Phi$ matrix}
\label{fig:connections_phi}
\end{center}

\vspace{-1mm} %layout

\subsection{IFOSMONDI-JFM's underlying non-linear problem}
\label{subsection:IFOSMONDI_JFM_s_underlying_non_linear_problem}

The IFOSMONDI-JFM co-simulation method \cite{Eguillon2021IfosmondiJFM} is based on the following principle: on each macro-step $[t^{[N]}, t^{[N+1]}[$, the inputs are defined so that their value and time-derivative at $t^{[N+1]}$ correspond to the connected outputs' value and time-derivative at $t^{[N+1]}$ too (hence, after the integration of the systems on this step). This will be called the \textit{coupling constraint} and is presented in \eqref{eq:IFOSMONDIJFM_coupling_constraint}.

\vspace{-3mm} %layout

\begin{equation}
\label{eq:IFOSMONDIJFM_coupling_constraint}
\renewcommand{\arraystretch}{1.7}
\left\{
	\begin{array}{lcl}
		\ub^{[N]}(t^{[N+1]}) & = & \Phi^T\ \tilde{\yb}^{[N+1]} \\
		\dspfrac{d\ub^{[N]}}{dt}(t^{[N+1]}) & = & \Phi^T\ \tilde{\dot{\yb}}^{[N+1]} \\
	\end{array}
\right.
\renewcommand{\arraystretch}{1.0}
\end{equation}

As the outputs and their derivatives are given, at $t^{[N+1]}$, by the functions $S_k^{[N]}$ and $\dot{S}_k^{[N]}$ on every system (for $k\in\Insys$), and as these functions require the expressions of all inputs, this problem is implicit. IFOSMONDI-JFM uses Newton-like methods to solve this implicit problem, see details in \cite{Eguillon2021IfosmondiJFM}. As each evaluation of the $S_k^{[N]}$ functions with different inputs require an integration of the systems, the method requires the rollback capability (see \ref{subsection:Capabilities}).

The inputs are defined so that they are $C^1$ on the whole co-simulation time domain, even at the communication times. In the middle of a co-simulation step, there is no problem: inputs are polynomials so they are $C^{\infty}$. However, for the inputs to be $C^1$ at the communication times, the method constraints their value and first-order derivative.

Therefore, when a communication time $t^{[N]}$ has been reached (with $N>0$), the value and time-derivative of every input are known at $t^{[N]}$ as, at this stage, the step $[t^{[N-1]}, t^{[N]}[$ has been integrated and validated (the method moved forward in time from $[t^{[N-1]}, t^{[N]}[$ to $[t^{[N]}, t^{[N+1]}[$). The inputs (on all systems) used at the valid and final integration on $[t^{[N-1]}, t^{[N]}[$ will be denoted $\ub_{\text{valid}}^{[N-1]}$. The $C^1$ smoothness at the beginning of $[t^{[N]}, t^{[N+1]}[$ (thus at $t^{[N]}$ on the right) is hence guaranteed by the constrains at $t^{[N]}$ in \eqref{eq:IFOSMONDIJFM_inputs_definition}.

To solve the coupling constraint \eqref{eq:IFOSMONDIJFM_coupling_constraint}, the Newton-like method in IFOSMONDI-JFM is used. The problem formulation is a zero finding of the function $\eta^{[N]}$ defined in \eqref{eq:IFOSMONDIJFM_etaN}

\vspace{-1mm} %layout

%layouted equation: several \! added
\begin{equation}
\label{eq:IFOSMONDIJFM_etaN}
\eta^{[N]}:
\left\{
	\begin{array}{lcl}
		\rnintot\times\rnintot & \!\!\rightarrow\!\! & \rnintot\times\rnintot \\
		\matDUbar{1.5}{\hat{\ub}^{[N+1]}}{\hat{\dot{\ub}}^{[N+1]}}
		& \!\!\mapsto\!\! &
		\matDUbar{1.5}{\hat{\ub}^{[N+1]}}{\hat{\dot{\ub}}^{[N+1]}}
			-
			\matDDbar{1.5}{\Phi^T}{0}{0}{\Phi^T}
			\begin{small}
			\left(
				\renewcommand{\arraystretch}{1.5}
				\begin{array}{c}
					S_1^{[N]}(\tilde{x}_1^{[N]}, \underbrace{u_1^{[N]}}) \\
					\vdots \\
					S_{\nsys}^{[N]}(\tilde{x}_{\nsys}^{[N]}, \underbrace{u_{\nsys}^{[N]}}) \\
					\hline
					\dot{S}_1^{[N]}(\tilde{x}_1^{[N]}, \underbrace{u_1^{[N]}}) \\
					\vdots \\
					\dot{S}_{\nsys}^{[N]}(\tilde{x}_{\nsys}^{[N]}, \underbrace{u_{\nsys}^{[N]}}) \\
				\end{array}
				\renewcommand{\arraystretch}{1.0}
			\right)
			\end{small}
			\hspace{2cm}\phantom{a}
		\vspace{-4mm}
		\\
		& & \multicolumn{1}{r}{\text{depend on }\begin{small}\matDUbar{1.5}{\hat{\ub}^{[N+1]}}{\hat{\dot{\ub}}^{[N+1]}}\end{small}}
	\end{array}
\right.
\end{equation}

\noindent where the zeros in the composite connection matrix denotes the null matrix $0_{\nintot\times\nouttot}$ of size $\nintot\times\nouttot$, where a vertical or horizontal bar represents vector or matrix concatenation, and where the $u_k^{[N]}$ time-dependent inputs for every $k\in\Insys$ are subparts of the total time-dependent input vector (as defined in \eqref{eq:ios_tot_fcts}). The latter is defined by the Hermite interpolation \eqref{eq:IFOSMONDIJFM_inputs_definition}.

\vspace{-3mm} %layout

\begin{equation}
\label{eq:IFOSMONDIJFM_inputs_definition}
\renewcommand{\arraystretch}{2.3}
\begin{array}{c|l|l}
	\multirow{2}{*}
	{$
		\ub^{[N]} \in \left(\mathbb{R}_3[t]\right)^{\nintot}
	$}
	&
	\ub^{[N]}(t^{[N]}) = \ub_{\text{valid}}^{[N-1]}(t^{[N]})
	&
	\dspfrac{d\ub^{[N]}}{dt}(t^{[N]}) = \dspfrac{d\ub_{\text{valid}}^{[N-1]}}{dt}(t^{[N]})
	\\
	&
	\ub^{[N]}(t^{[N+1]}) = \hat{\ub}^{[N+1]}
	&
	\dspfrac{d\ub^{[N]}}{dt}(t^{[N+1]}) = \hat{\dot{\ub}}^{[N+1]}
\end{array}
\renewcommand{\arraystretch}{1.0}
\end{equation}

\vspace{-1mm} %layout

\noindent where $\hat{\ub}^{[N+1]}$ and $\hat{\dot{\ub}}^{[N+1]}$ are the variables of $\eta^{[N]}$ (see \eqref{eq:IFOSMONDIJFM_etaN}).

Hence, once the root of $\eta^{[N]}$ is found by the Newton-like method, it is used to define $\ub_{\text{valid}}^{[N]}$ with the Hermite interpolation \eqref{eq:IFOSMONDIJFM_inputs_definition}. This $\ub_{\text{valid}}^{[N]}$ satisfies the coupling constraint \eqref{eq:IFOSMONDIJFM_coupling_constraint}, and the validated inputs are $C^1$ in $t^{[N]}$.

More details about the underlying computations (namely for the first macro-step and the starting point of the zero-finding of $\eta^{[N]}$ on each new macro-step) are given on reference \cite{Eguillon2021IfosmondiJFM} which focuses on the IFOSMONDI-JFM method.

\section{MISSILES method}
\label{section:MISSILES_method}

Based on the formalism introduced in \ref{section:Formalism_and_notations} and on the willingness to get benefit from the benefits of the IFOSMONDI-JFM method without the rollback capability requirement, we define in this section the MISSILES co-simulation method, standing for \textit{Mock Iteration for Solving Smooth Interfaces with Linear Estimations of Systems}.

\subsection{General idea}
\label{subsection:general_idea}

Of course, a first idea could be to simply replace the non-rollback capable systems by a COSTARICA (it is indeed the aim of the COSTARICA process). However, in case every system is replaced by a COSTARICA, the relationship between the inputs at a given macro-step and the outputs at the end of this macro-step is known (cf. subsection \ref{subsection:COSTARICA_estimator}). We can therefore replace the occurences of $S_k^{[N]}$ and $\dot{S}_k^{[N]}$ for all systems $k\in\Insys$ in the expression of $\eta^{[N]}$ (the function whose root is searched, in \eqref{eq:IFOSMONDIJFM_etaN}) and use the resulting function expression to directly find its root. This modified version of $\eta^{[N]}$ will be denoted by $\eta_{\text{MISSILES}}^{[N]}$:

\vspace{-6mm} %layout

%layouted equation: several \! added
\begin{equation}
\label{eq:eta_MISSILES_N}
\eta_{\text{MISSILES}}^{[N]}\!:
\left\{\!
	\begin{array}{lcl}
		\rnintot\times\rnintot & \!\!\!\!\!\!\rightarrow\!\!\!\! & \rnintot\times\rnintot \\
		\matDUbar{1.5}{\hat{\ub}^{[N+1]}}{\hat{\dot{\ub}}^{[N+1]}}
		& \!\!\!\!\!\!\mapsto\!\!\!\! &
		\matDUbar{1.5}{\hat{\ub}^{[N+1]}}{\hat{\dot{\ub}}^{[N+1]}}
			-
			\matDDbar{1.5}{\Phi^T}{0}{0}{\Phi^T}
			\begin{small}
			\left(
				\renewcommand{\arraystretch}{1.5}
				\begin{array}{c}
					\GV_1^{[N]}\check{\Xid}_1^{[N]}+\PV_1^{[N]}\tilde{x}_1^{[N]}+\RV_1^{[N]}\tfC_1^{[N]} + \hyC_1^{[N+1]} \\
					\vdots \\
					\GV_{\nsys}^{[N]}\check{\Xid}_{\nsys}^{[N]}+\PV_{\nsys}^{[N]}\tilde{x}_{\nsys}^{[N]}+\RV_{\nsys}^{[N]}\tfC_{\nsys}^{[N]} + \hyC_{\nsys}^{[N+1]} \\
					\hline
					\GD_1^{[N]}\check{\Xid}_1^{[N]}+\PD_1^{[N]}\tilde{x}_1^{[N]}+\RD_1^{[N]}\tfC_1^{[N]} + \hdyC_1^{[N+1]} \\
					\vdots \\
					\GD_{\nsys}^{[N]}\check{\Xid}_{\nsys}^{[N]}+\PD_{\nsys}^{[N]}\tilde{x}_{\nsys}^{[N]}+\RD_{\nsys}^{[N]}\tfC_{\nsys}^{[N]} + \hdyC_{\nsys}^{[N+1]} \\
				\end{array}
				\renewcommand{\arraystretch}{1.0}
			\right)
			\end{small}
	\end{array}
\right.
\end{equation}

\vspace{-4mm} %layout

\noindent where the matrices $\check{\Xid}_1^{[N]}$, ..., $\check{\Xid}_{\nsys}^{[N]}$ introduced in \eqref{eq:COSTARICA_yL_is_GXi_plus_Pxi} depend on \begin{small}$\matDUbar{1.5}{\hat{\ub}^{[N+1]}}{\hat{\dot{\ub}}^{[N+1]}}$\end{small} as the latter is used to calibrate the polynomial inputs $\ub^{[N]}$ (see \eqref{eq:IFOSMONDIJFM_inputs_definition}) and as the $\check{\Xid}_1^{[N]}$, ..., $\check{\Xid}_{\nsys}^{[N]}$ matrices contain the coefficients of the time-shifted version of $\ub^{[N]}$ (see subsection \ref{subsection:time_shift} for details about time-shift). This relationship will be detailed further in \ref{subsection:global_linear_problem}.

MISSILES removes the iterative part of the IFOSMONDI-JFM method and replaces it by a single resolution (detailed further in this paper) to satisfy the coupling constraint. However, this resolution is based on the COSTARICA estimators on each system ("\textit{Linear Estimations of Systems}" in MISSILES' name come from here). Therefore, this is not really an iteration on the systems, but a single fake iteration as the systems are not integrated at this stage ("\textit{Mock Iteration}" in MISSILES' name come from here).

Once this resolution is done (root finding of $\eta_{\text{MISSILES}}^{[N]}$), the input values and time-derivatives satisfying the coupling constraint (on the COSTARICA surrogates) are use together with the input values and time-derivatives at the end of the previous macro-step (\textit{i.e.} beginning of the current macro-step) to directly define the validated input expressions $\ub_{\text{valid}}^{[N]}$. The systems can then be integrated on the macro-step $[t^{[N]}, t^{[N+1]}[$ parallely with their respective inputs. The whole process is shown on figure \ref{fig:MISSILES_method}.

\vspace{-2mm} %layout

\begin{center}
\includegraphics[scale=0.29]{\figuresdir/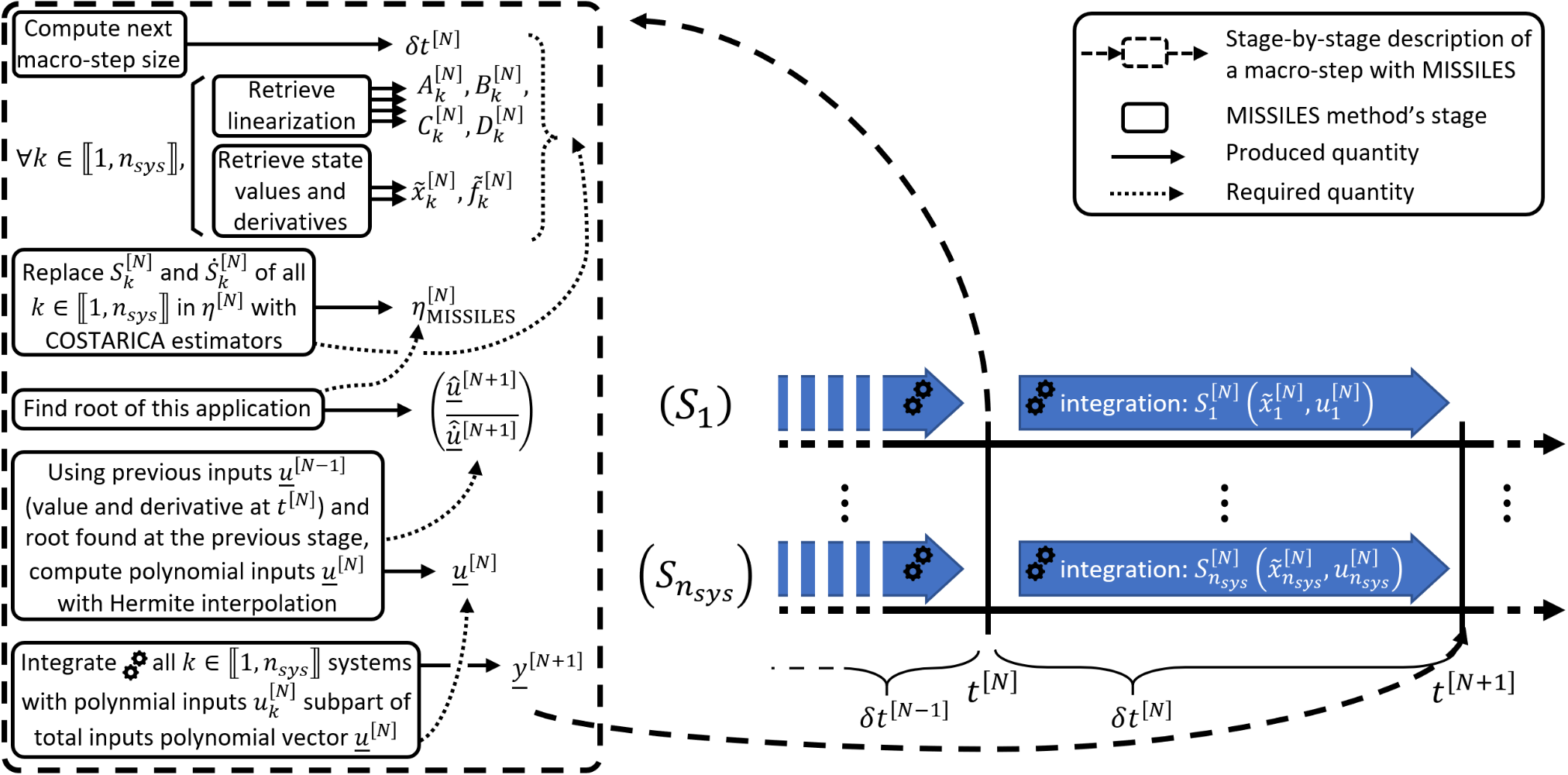}
\vspace{-2mm} %layout
\captionof{figure}{Schematic view of MISSILES co-simulation method}
\label{fig:MISSILES_method}
\end{center}

Three main characteristics of the MISSILES method can be noticed so far:

\vspace{-2mm} %layout

\begin{itemize}
\setlength{\itemsep}{-0.4em} %layout
\ite \underline{Synchronicity:} The communication times $(t^{[N]})_{N\in[\![0, \Nmax]\!]}$ are the same on every system. The method is said to be \textit{synchronous}, as well as the IFOSMONDI-JFM method namely (in opposition to \textit{asynchronous} co-simulation methods such as \cite{Muller2016} or \cite{Eguillon2022F3ornits}).
\ite \underline{Explicit nature:} The polynomial inputs of all systems are predicted on $[t^{[N]}, t^{[N+1]}[$ when when time $t^{[N]}$ is reached. This prediction is then used in every system for their single (first and last on this very macro-step) integration on the current macro-step. The rollback capability is therefore not required.
\ite \underline{Parallelizability:} The stage where the systems have to be integrated for real can be done in parallel in the sense that the inputs for each system can be computed using only data known when all systems have been simulated until time $t^{[N]}$, contrary to Gauss-Seidel-based co-simulation methods \cite{Burrage1995} or others like \cite{Benedikt2019TriggerSeq}.
\end{itemize}

\vspace{-2mm} %layout

From the required capabilities point of view, the rollback, mandatory in the case of IFOSMONDI-JFM, is no more required on MISSILES. However, the linearization and the state variables retrievals are now required due to the use of the COSTARICA estimators, despite these capabilities were not required in the case of the IFOSMONDI-JFM method. However, as the rollback is way scarser than the ability to retrieve the linearization and the state variables in practice, we can reasonably stand that the MISSILES co-simulation algorithm is more usable in an industrial context.

Table \ref{tab:capabilities_MISSILES} sums up the required available interactions in the co-simulation systems. Each interaction is referred to by its designation in table \ref{tab:capabilities}.

\vspace{-5mm} %layout

%layouted table: several \hspace{-3mm} added, + \slwordwrap environment to put text in 'small' + paddings rearrangements
\begin{center}
\captionof{table}{Interaction with the co-simulation systems for the MISSILES method and justification}
\label{tab:capabilities_MISSILES}
\renewcommand{\arraystretch}{1.5}
\begin{tabular}{|c|c|l|}
	\hline
	\wordwrap{1.2}{Interaction \\name} & \wordwrap{1.2}{Is it \\required?} & \hfil Why? (justification) \hfil \\
	\hline
	\hline
	Provide inputs & Yes & \hspace{-3mm}\slwordwrap{small}{1.2}{Basic interaction required for all co-simulation methods}\hspace{-3mm} \\
	\hline
	\wordwrap{1.2}{Provide time-\\dependent inputs} & Yes & \hspace{-3mm}\slwordwrap{small}{1.2}{As in IFOSMONDI-JFM, the inputs must be $C^1$ at the communication \\times, satisfy the coupling constraint at $t^{[N+1]}$ on values and time-\\derivatives, and at $t^{[N]}$ (\textit{de facto}, due to their $C^1$ character). Constant \\functions cannot do so on non-constant signals.}\hspace{-3mm} \\
	\hline
	Do a step & Yes & \hspace{-3mm}\slwordwrap{small}{1.2}{Basic interaction required for all co-simulation methods}\hspace{-3mm} \\
	\hline
	Retrieve outputs & Yes & \hspace{-3mm}\slwordwrap{small}{1.2}{Basic interaction required for all co-simulation methods}\hspace{-3mm} \\
	\hline
	\wordwrap{1.2}{Retrieve outputs \\time-derivatives} & No & \hspace{-3mm}\slwordwrap{small}{1.2}{Despite IFOSMONDI-JFM needs to know the outputs time-derivatives \\to evaluate the coupling constraint, the COSTARICA estimators can \\estimate them based on $\GD_k^{[N]}$ and $\PD_k^{[N]}$ on all systems.\padDown}\hspace{-3mm} \\
	\hline
	\hspace{-4mm}
	\wordwrap{1.2}{Retrieve state values \\and derivatives}
	\hspace{-4mm} & Yes\vphantom{\wordwrap{1.2}{\PadUp a \PadDown}} & \multirow{2}{*}{\hspace{-3mm}\slwordwrap{small}{1.2}{\padUp The substitution of the $S_k^{[N]}$ and $\dot{S}_k^{[N]}$ for all systems $k\in\Insys$ by the \\COSTARICA estimators in $\eta_{\text{MISSILES}}^{[N]}$ makes MISSILES require the \\states and the $A^{[N]}$, $B^{[N]}$, $C^{[N]}$ and $D^{[N]}$ matrices to compute these estimators.}\hspace{-3mm}} \\
	\cline{1-2}
	Retrieve linearization & Yes\vphantom{\wordwrap{1.2}{\PadUp a \PadDown}} & \\
	\hline
	Rollback & No & \hspace{-3mm}\slwordwrap{small}{1.2}{\PadUp Contrary to evaluations of $\eta^{[N]}$, the evaluations of $\eta_{\text{MISSILES}}^{[N]}$ do not \\require to integrate the systems. The root finding process hence does \\not require to roll back the systems to re-evaluate the function.}\hspace{-3mm} \\
	\hline
\end{tabular}
\renewcommand{\arraystretch}{1.0}
\end{center}

\vspace{-2mm} %layout

The remaining question is: \textbf{How to find the root of $\eta_{\text{MISSILES}}^{[N]}$?} Indeed, on figure \ref{fig:MISSILES_method}, the stage called "find root of this application" hasn't been discussed so far. Following subsection is dedicated to the assembly of a linear problem which solution is the root of $\eta_{\text{MISSILES}}^{[N]}$.

\vspace{-2mm} %layout

\subsection{Global linear problem}
\label{subsection:global_linear_problem}

Let's gather the elements introduced previously in section \ref{section:Formalism_and_notations} in order to find the root of $\eta_{\text{MISSILES}}^{[N]}$. Let's consider the solution \begin{small}$\matDUbar{1.5}{\hat{\ub}^{[N+1]}}{\hat{\dot{\ub}}^{[N+1]}}$\end{small} to the root finding problem. By definition of $\eta_{\text{MISSILES}}^{[N]}$ in \eqref{eq:eta_MISSILES_N}, we have:

%layouted figure: several \vspace{-1mm} added
\begin{equation}
\label{eq:eta_MISSILES_N_solution}
\eta_{\text{MISSILES}}^{[N]}\matDUbar{1.5}{\hat{\ub}^{[N+1]}}{\hat{\dot{\ub}}^{[N+1]}}
=
0
\hsd
\Leftrightarrow
\hsd
\matDUbar{1.5}{\hat{\ub}^{[N+1]}}{\hat{\dot{\ub}}^{[N+1]}}
=
\matDDbar{1.5}{\Phi^T}{0}{0}{\Phi^T}
\begin{small}
\left(
	\renewcommand{\arraystretch}{1.5}
	\begin{array}{c}
		\GV_1^{[N]}\check{\Xid}_1^{[N]}+\PV_1^{[N]}\tilde{x}_1^{[N]}+\RV_1^{[N]}\tfC_1^{[N]} + \hyC_1^{[N+1]} \vspace{-1mm}\\
		\vdots \\
		\GV_{\nsys}^{[N]}\check{\Xid}_{\nsys}^{[N]}+\PV_{\nsys}^{[N]}\tilde{x}_{\nsys}^{[N]}+\RV_{\nsys}^{[N]}\tfC_{\nsys}^{[N]} + \hyC_{\nsys}^{[N+1]} \\
		\hline
		\GD_1^{[N]}\check{\Xid}_1^{[N]}+\PD_1^{[N]}\tilde{x}_1^{[N]}+\RD_1^{[N]}\tfC_1^{[N]} + \hdyC_1^{[N+1]} \vspace{-1mm}\\
		\vdots \\
		\GD_{\nsys}^{[N]}\check{\Xid}_{\nsys}^{[N]}+\PD_{\nsys}^{[N]}\tilde{x}_{\nsys}^{[N]}+\RD_{\nsys}^{[N]}\tfC_{\nsys}^{[N]} + \hdyC_{\nsys}^{[N+1]} \\
	\end{array}
	\renewcommand{\arraystretch}{1.0}
\right)
\end{small}
\end{equation}

\vspace{-4mm} %layout

Let's now consider the outputs (and their time-derivatives) estimations at $t^{[N+1]}$, denoted by \begin{small}$\matDUbar{1.5}{\hat{\yb}^{[N+1]}}{\hat{\dot{\yb}}^{[N+1]}}$\end{small} which, when dispatched according to the modular model's topology, generate the solution to \eqref{eq:eta_MISSILES_N_solution}. Due to the structure of the $\Phi$ matrix \eqref{eq:Phi} presented in \ref{subsection:connecting_systems_into_a_modular_model} (either permutation or not, by with a single $1$ on each and every column), to a given \begin{small}$\matDUbar{1.5}{\hat{\ub}^{[N+1]}}{\hat{\dot{\ub}}^{[N+1]}}$\end{small} corresponds a single\footnote{The reciprocal is not always true: in case an output is connected to several inputs, a row of $\Phi$ has several $1$ coefficients, as in figure \ref{fig:connections_phi}. In this case, if the inputs connected to the same output have different values, no output vector corresponds to the input vector.} \begin{small}$\matDUbar{1.5}{\hat{\yb}^{[N+1]}}{\hat{\dot{\yb}}^{[N+1]}}$\end{small}. We can thus calculate the solution on the outputs instead of the inputs. The solution to \eqref{eq:eta_MISSILES_N_solution} will then be retrievable with \eqref{eq:MISSILES_solution_dispatching}.

\vspace{-4mm} %layout

\begin{equation}
\label{eq:MISSILES_solution_dispatching}
\matDUbar{1.5}{\hat{\ub}^{[N+1]}}{\hat{\dot{\ub}}^{[N+1]}}
=
\matDDbar{1.5}{\Phi^T}{0}{0}{\Phi^T}
\matDUbar{1.5}{\hat{\yb}^{[N+1]}}{\hat{\dot{\yb}}^{[N+1]}}
\end{equation}

\vspace{-2mm} %layout

Therefore, to generate the solution to \eqref{eq:eta_MISSILES_N_solution}, the outputs solutions must write:

\vspace{-4mm} %layout

%layouted figure: several \vspace{-1mm} added
\begin{equation}
\label{eq:eta_MISSILES_N_solution_outputs}
\matDUbar{1.5}{\hat{\yb}^{[N+1]}}{\hat{\dot{\yb}}^{[N+1]}}
=
\begin{small}
\left(
	\renewcommand{\arraystretch}{1.5}
	\begin{array}{c}
		\GV_1^{[N]}\check{\Xid}_1^{[N]}+\PV_1^{[N]}\tilde{x}_1^{[N]}+\RV_1^{[N]}\tfC_1^{[N]} + \hyC_1^{[N+1]} \vspace{-1mm}\\
		\vdots \\
		\GV_{\nsys}^{[N]}\check{\Xid}_{\nsys}^{[N]}+\PV_{\nsys}^{[N]}\tilde{x}_{\nsys}^{[N]}+\RV_{\nsys}^{[N]}\tfC_{\nsys}^{[N]} + \hyC_{\nsys}^{[N+1]} \\
		\hline
		\GD_1^{[N]}\check{\Xid}_1^{[N]}+\PD_1^{[N]}\tilde{x}_1^{[N]}+\RD_1^{[N]}\tfC_1^{[N]} + \hdyC_1^{[N+1]} \vspace{-1mm} \\
		\vdots \\
		\GD_{\nsys}^{[N]}\check{\Xid}_{\nsys}^{[N]}+\PD_{\nsys}^{[N]}\tilde{x}_{\nsys}^{[N]}+\RD_{\nsys}^{[N]}\tfC_{\nsys}^{[N]} + \hdyC_{\nsys}^{[N+1]} \\
	\end{array}
	\renewcommand{\arraystretch}{1.0}
\right)
\end{small}
\end{equation}

\vspace{-2.5mm} %layout

The solution of problem \eqref{eq:eta_MISSILES_N_solution_outputs} leads, with dispatching \eqref{eq:MISSILES_solution_dispatching}, to the solution of \eqref{eq:eta_MISSILES_N_solution}. Let's detail the five stage to get, from the $\check{\Xid}_1^{[N]}$, ..., $\check{\Xid}_{\nsys}^{[N]}$ matrices in \eqref{eq:eta_MISSILES_N_solution_outputs} (introduced in \eqref{eq:COSTARICA_yL_is_GXi_plus_Pxi} and detailed in \ref{subsection:time_shift}), the outputs \begin{small}$\matDUbar{1.5}{\hat{\yb}^{[N+1]}}{\hat{\dot{\yb}}^{[N+1]}}$\end{small}:

\vspace{-3.5mm} %layout

\begin{enumerate}
\setlength{\itemsep}{0em} %layout
\item In \eqref{eq:eta_MISSILES_N_solution_outputs}, the $\check{\Xid}_1^{[N]}$, ..., $\check{\Xid}_{\nsys}^{[N]}$ matrices contain the coefficients of the time-shifted version of the time-dependent inputs, as mentioned earlier. For all system $k\in\Insys$, this time-shift corresponds to the tensor-matrix product $\check{\Xid}_k^{[N]} = \mathcal{C}_k^{[N]} \Xid_k^{[N]}$ introduced in \eqref{eq:Xi_check_is_CkN_times_Xi}.
\item The $\Xid_k^{[N]}$ matrix mentioned in step 1 contains the coefficients of the polynomial inputs $u_k^{[N]}$ (see \eqref{eq:Xi_and_ukN_coefficients}).
\item The polynomial inputs $u_k^{[N]}$ mentioned in step 2 are calibrated with the Hermite interpolation \eqref{eq:IFOSMONDIJFM_inputs_definition} (the latter acts on the total input vector, yet we can do it system by system as subparts of the vector $\ub^{[N]}$, see \eqref{eq:ios_tot_fcts}).
\item The coefficients of the inputs can be expressed linearly with respect to the input constraints at $t^{[N+1]}$ of the Hermite interpolation mentioned in step 3. This linear expression has been introduced in subsection \ref{subsection:hermite_interpolation}. The matrix and vector or this linear expression only depend on the times and the input constraints at $t^{[N]}$, known and independent of \begin{small}$\matDUbar{1.5}{\hat{\yb}^{[N+1]}}{\hat{\dot{\yb}}^{[N+1]}}$\end{small} due to the $C^1$ condition.
\item Finally, the input constraints at $t^{[N+1]}$ mentioned in step 4 can be obtained from the output constraints at the same time, that is to say \begin{small}$\matDUbar{1.5}{\hat{\yb}^{[N+1]}}{\hat{\dot{\yb}}^{[N+1]}}$\end{small}, using the dispatching relationship \eqref{eq:MISSILES_solution_dispatching}.
\end{enumerate}

The application of these steps to equation \eqref{eq:eta_MISSILES_N_solution_outputs} gives:

\vspace{-4mm} %layout

\begin{equation}
\label{eq:eta_MISSILES_N_solution_outputs_developed}
\matDUbar{1.5}{\hat{\yb}^{[N+1]}}{\hat{\dot{\yb}}^{[N+1]}}
=
\matDUbar{1.5}{\GbV^{[N]}}{\GbD^{[N]}}
\bb{\mathcal{C}}^{[N]}
\left(
	\bb{\mathcal{A}}^{[N]}
	\matDDbar{1.5}{\Phi^T}{0}{0}{\Phi^T}		
	\matDUbar{1.5}{\hat{\yb}^{[N+1]}}{\hat{\dot{\yb}}^{[N+1]}}
	+
	\bb{\mathcal{B}}^{[N]}
\right)
+
\matDUbar{1.5}{\PbV^{[N]}}{\PbD^{[N]}}
\bb{\tilde{x}}^{[N]}
+
\matDUbar{1.5}{\RbV^{[N]}}{\RbD^{[N]}}
\bb{\tfC}^{[N]}
+
\matDUbar{1.5}{\bb{\hat{y}_C}^{[N+1]}}{\bb{\hat{\dot{y}}_C}^{[N+1]}}
\end{equation}

\vspace{-1mm} %layout

\noindent where the underline tensors and matrices are composition of previously introduced by-system quantities. These global operators are described below using by-system operators introduced in section \ref{section:Formalism_and_notations}. Please note that, despite the maximum input polynomial degree is $n=3$ (due to the Hermite interpolation \eqref{eq:IFOSMONDIJFM_inputs_definition}), is it still written as $n$ below for the sake of genericity. Let's recall that the degree is the only dimension in the tensors that starts by zero so that the $p$\up{th} element correspond to the monomial of degree $p$.

\begin{small}$\matDUbar{1.5}{\GbV^{[N]}}{\GbD^{[N]}}$\end{small} is a tensor of order $3$ and of size $2\ \nouttot\times\nintot\times (n+1)$. It represents the action of all inputs in the COSTARICA estimators, and is a concatenation of the two tensors $\GbV^{[N]}$ and $\GbD^{[N]}$, both of order $3$ and of size $\nouttot\times\nintot\times (n+1)$.

\vspace{-7mm} %layout

\begin{equation}
\label{eq:GlobOp_Gb}
\forall (\ib, \jb, p)\in[\![1, 2\ \nouttot]\!] \times \Inintot \times \Zn,\
\matDUbar{1.5}{\GbV^{[N]}}{\GbD^{[N]}}_{\ib, \jb, p}
\!\!\!\!\!\! %layout
=
\left\{
	\begin{array}{ll}
		\left(
			\GbV^{[N]}
		\right)_{\ib, \jb, p}
		&
		\hspace{-3mm} %layout
		\text{ if } \ib\in\Inouttot
		\\
		\left(
			\GbD^{[N]}
		\right)_{(\ib-\nouttot), \jb, p}
		&
		\hspace{-3mm} %layout
		\text{ if } \ib\in[\![\nouttot+1, 2\ \nouttot]\!]
	\end{array}
\right.
\end{equation}

\vspace{-4mm} %layout

\begin{equation}
\label{eq:GlobOp_GbV_and_Gbd}
\begin{array}{r}
	\multicolumn{1}{l}
	{
		\forall k\in\Insys,\
		\forall l\in\Insys,\
		\forall i\in\Inoutk{k},\
		\forall j\in\Inink{l},\
		\forall p\in\Zn,\
		\hspace{3cm}
		\phantom{a}
	}
	\\
	\left(\GbV^{[N]}\right)_{
		(i + \sum_{\kappa=1}^{k-1} \noutk{\kappa}),
		(j + \sum_{\lambda=1}^{l-1} \nink{\lambda}),
		p
	} = \delta_{k, l} \cdot \left(\GV_k^{[N]}\right)_{i, j, p}
	\\
	\text{and }
	\left(\GbD^{[N]}\right)_{
		(i + \sum_{\kappa=1}^{k-1} \noutk{\kappa}),
		(j + \sum_{\lambda=1}^{l-1} \nink{\lambda}),
		p
	} = \delta_{k, l} \cdot \left(\GD_k^{[N]}\right)_{i, j, p}
	\\
\end{array}
\end{equation}

\vspace{-5mm} %layout

\begin{small}$\matDUbar{1.5}{\PbV^{[N]}}{\PbD^{[N]}}$\end{small} is a matrix of size $2\ \nouttot\times\nsttot$. It represents the effect of the initial states on the current macro-step in the COSTARICA estimators, and is a concatenation of the two matrices $\PbV^{[N]}$ and $\PbD^{[N]}$, both of size $\nouttot\times\nsttot$.

\vspace{-6mm} %layout

\begin{equation}
\label{eq:GlobOp_Pb}
\forall (\ib, \sigmab)\in[\![1, 2\ \nouttot]\!] \times \Insttot,\
\matDUbar{1.5}{\PbV^{[N]}}{\PbD^{[N]}}_{\ib, \sigmab}
=
\left\{
	\begin{array}{ll}
		\left(
			\PbV^{[N]}
		\right)_{\ib, \sigmab}
		&
		\text{ if } \ib\in\Inouttot
		\\
		\left(
			\PbD^{[N]}
		\right)_{(\ib-\nouttot), \sigmab}
		&
		\text{ if } \ib\in[\![\nouttot+1, 2\ \nouttot]\!]
	\end{array}
\right.
\end{equation}

\vspace{-2mm} %layout

\begin{equation}
\label{eq:GlobOp_PbV_and_PbD}
\begin{array}{r}
	\multicolumn{1}{l}
	{
		\forall k\in\Insys,\
		\forall l\in\Insys,\
		\forall i\in\Inoutk{k},\
		\forall \sigma\in\Instk{l},\
		\hspace{3cm}
		\phantom{a}
	}
	\\
	\left(\PbV^{[N]}\right)_{
		(i + \sum_{\kappa=1}^{k-1} \noutk{\kappa}),
		(\sigma + \sum_{\lambda=1}^{l-1} \nstk{\lambda})
	} = \delta_{k, l} \cdot \left(\PV_k^{[N]}\right)_{i, \sigma}
	\\
	\text{and }
	\left(\PbD^{[N]}\right)_{
		(i + \sum_{\kappa=1}^{k-1} \noutk{\kappa}),
		(\sigma + \sum_{\lambda=1}^{l-1} \nstk{\lambda})
	} = \delta_{k, l} \cdot \left(\PD_k^{[N]}\right)_{i, \sigma}
	\\
\end{array}
\end{equation}

\vspace{-5mm} %layout

\begin{small}$\matDUbar{1.5}{\RbV^{[N]}}{\RbD^{[N]}}$\end{small} is a matrix of size $2\ \nouttot\times\nsttot$. It represents the effect of the difference $\left(\tfC_k\right)_{k\in\Insys}$ between the linearization and the state-space representation of the systems (see \cite{Eguillon2022Costarica}). It is a concatenation of the two matrices $\RbV^{[N]}$ and $\RbD^{[N]}$, both of size $\nouttot\times\nsttot$.

\vspace{-6mm} %layout

\begin{equation}
\label{eq:GlobOp_Rb}
\forall (\ib, \sigmab)\in[\![1, 2\ \nouttot]\!] \times \Insttot,\
\matDUbar{1.5}{\RbV^{[N]}}{\RbD^{[N]}}_{\ib, \sigmab}
=
\left\{
	\begin{array}{ll}
		\left(
			\RbV^{[N]}
		\right)_{\ib, \sigmab}
		&
		\text{ if } \ib\in\Inouttot
		\\
		\left(
			\RbD^{[N]}
		\right)_{(\ib-\nouttot), \sigmab}
		&
		\text{ if } \ib\in[\![\nouttot+1, 2\ \nouttot]\!]
	\end{array}
\right.
\end{equation}

\vspace{-2mm} %layout

\begin{equation}
\label{eq:GlobOp_RbV_and_RbD}
\begin{array}{r}
	\multicolumn{1}{l}
	{
		\forall k\in\Insys,\
		\forall l\in\Insys,\
		\forall i\in\Inoutk{k},\
		\forall \sigma\in\Instk{l},\
		\hspace{3cm}
		\phantom{a}
	}
	\\
	\left(\RbV^{[N]}\right)_{
		(i + \sum_{\kappa=1}^{k-1} \noutk{\kappa}),
		(\sigma + \sum_{\lambda=1}^{l-1} \nstk{\lambda})
	} = \delta_{k, l} \cdot \left(\RV_k^{[N]}\right)_{i, \sigma}
	\\
	\text{and }
	\left(\RbD^{[N]}\right)_{
		(i + \sum_{\kappa=1}^{k-1} \noutk{\kappa}),
		(\sigma + \sum_{\lambda=1}^{l-1} \nstk{\lambda})
	} = \delta_{k, l} \cdot \left(\RD_k^{[N]}\right)_{i, \sigma}
	\\
\end{array}
\end{equation}

\vspace{-2mm} %layout

$\bb{\tfC}^{[N]}$ is a column vector of size $\nsttot$. It corresponds to a concatenation of the $\tfC_k$ parts in the COSTARICA estimators introduced in \ref{subsection:COSTARICA_estimator}, equation \eqref{eq:COSTARICA_tfC}.

\vspace{-2mm} %layout

\begin{equation}
\label{eq:GlobOp_tfC}
\renewcommand{\arraystretch}{1.9}
\begin{array}{lccccr}
	\bb{\tfC}^{[N]} =
	\Big(
		&
		\tfC_{1, 1}^{[N]}\ ,\ ...\ ,\ \tfC_{1, \nstk{1}}^{[N]}\ ,
		\hsd & \hsd
		\tfC_{2, 1}^{[N]}\ ,\ ...\ ,\ \tfC_{2, \nstk{2}}^{[N]}\ ,
		\hsd & \hsd
		...\ ,
		\hsd & \hsd
		\tfC_{\nsys, 1}^{[N]}\ ,\ ...\ ,\ \tfC_{\nsys, \nstk{\nsys}}^{[N]}
		&
	\Big)^T
\end{array}
\renewcommand{\arraystretch}{1.0}
\end{equation}

$\bb{\mathcal{A}}^{[N]}$ is a tensor of order $3$ and of size $\nintot\times (n+1)\times 2 \nintot$, and $\bb{\mathcal{B}}^{[N]}$ is a matrix of size $\nintot\times (n+1)$. They represent the Hermite interpolation, transforming the constraints upon inputs values and derivatives at the end of the current macro-step into the coefficients of the polynomial inputs on this step. They are compositions of $\AVelem^{[N]}$ and $\ADelem^{[N]}$ elements and $\Belem^{[N]}$ evaluations. These quantities are defined in \eqref{eq:GlobOp_AVelem_ADelem_Belem_on_tN_tNp1} from concrete applications on $[t^{[N]}, t^{[N+1]}[$ of the quantities introduced in \eqref{eq:Hermite_linexp_general} \eqref{eq:Hermite_AVelem_ADelem_Belem} in subsection \ref{subsection:hermite_interpolation}.

\vspace{-5mm} %layout

\begin{equation}
\label{eq:GlobOp_AVelem_ADelem_Belem_on_tN_tNp1}
\text{With }
\left\{
	\begin{array}{ll}
		t_1 & = t^{[N]} \\
		t_2 & = t^{[N+1]} \\
	\end{array}
\right\}
\text{ we define }
\left\{
	\begin{array}{l}
		\left\{
			\begin{array}{ll}
				\AVelem^{[N]} & \overset{\Delta}{=} \AVelem \\
				\ADelem^{[N]} & \overset{\Delta}{=} \ADelem \\
			\end{array}
		\right\}
		\text{ from definition \eqref{eq:Hermite_AVelem_ADelem_Belem}}
		\\
		\Belem^{[N]}:
		\begin{array}{lcl}
			(l, j)
			& \mapsto &
			\Belem
			\text{ as in \eqref{eq:Hermite_AVelem_ADelem_Belem} with }
			\left\{
				\begin{array}{ll}
					v_1 & = u_{l, j}^{[N-1]}(t^{[N]}) \\
					\dot{v}_1 & = \dspfrac{du_{l, j}^{[N-1]}}{dt}(t^{[N]}) \\
				\end{array}
			\right.
			\\
			& & 
			\forall (l, j) \in \left\{(l, j) |	l\in\Insys \text{ and } j\in\Inink{l}\right\}
		\end{array}
	\end{array}
\right.
\end{equation}

\vspace{-1mm} %layout

With \eqref{eq:GlobOp_AVelem_ADelem_Belem_on_tN_tNp1}, we apply the generic problem \eqref{eq:Hermite_2points_generic_properties} to the case \eqref{eq:IFOSMONDIJFM_inputs_definition}, where constraints $v_2$ and $\dot{v}_2$ in \eqref{eq:Hermite_2points_generic_properties} represent the solution of \eqref{eq:eta_MISSILES_N_solution} we are looking for. Global composite tensor $\bb{\mathcal{A}}^{[N]}$ is then defined in \eqref{eq:GlobOp_Acal} and global matrix $\bb{\mathcal{B}}^{[N]}$ in \eqref{eq:GlobOp_Bcal}.

\vspace{-3mm} %layout

\begin{equation}
\label{eq:GlobOp_Acal}
\begin{array}{r}
	\multicolumn{1}{l}
	{
		\forall l_1\in\Insys,\
		\forall j_1\in\Inink{l_1}\,
		\forall p\in\Zn,\
		\forall l_2\in\Insys,\
		\forall j_2\in\Inink{l_2},\
		\hspace{3cm}
		\phantom{a}
	}
	\\
	\begin{array}{rll}
		&
		\left(
			\bb{\mathcal{A}}^{[N]}
		\right)_{
			(j_1 + \sum_{\lambda=1}^{l_1-1} \nink{\lambda}),
			p,
			(j_2 + \sum_{\lambda=1}^{l_2-1} \nink{\lambda})
		} & =
		\left\{
			\begin{array}{ll}
				\left(\AVelem\right)_p & \text{ if } (l_1, j_1) = (l_2, j_2) \\
				0 & \text{ otherwise}
			\end{array}
		\right.
		\\
		\text{and }
		&
		\left(
			\bb{\mathcal{A}}^{[N]}
		\right)_{
			(j_1 + \sum_{\lambda=1}^{l_1-1} \nink{\lambda}),
			p,
			\nintot + (j_2 + \sum_{\lambda=1}^{l_2-1} \nink{\lambda})
		} & =
		\left\{
			\begin{array}{ll}
				\left(\ADelem\right)_p & \text{ if } (l_1, j_1) = (l_2, j_2) \\
				0 & \text{ otherwise}
			\end{array}
		\right.
	\end{array}
\end{array}
\end{equation}

\vspace{-2mm} %layout

\begin{equation}
\label{eq:GlobOp_Bcal}
\forall l\in\Insys,\
\forall j\in\Inink{l}\,
\forall p\in\Zn,\
\left(
	\bb{\mathcal{B}}^{[N]}
\right)_{
	(j + \sum_{\lambda=1}^{l-1} \nink{\lambda}),
	p
} =	\left(\Belem\right)_p
\end{equation}

\vspace{-1mm} %layout

$\bb{\mathcal{C}}^{[N]}$ is a tensor of order $4$ and of size $\nintot\times (n+1)\times\nintot\times (n+1)$. It represents the time-shift of the coefficients of all polynomial inputs from $[t^{[N]}, t^{[N+1]}[$ to $[0, \dt^{[N]}[$ as described in \ref{subsection:time_shift}.

\vspace{-3mm} %layout

%layouted figure: \padDown added and several \hspace{-2mm}
\begin{equation}
\label{eq:GlobOp_Ccal}
\begin{array}{r}
	\multicolumn{1}{l}
	{
		\forall l_1\in\Insys,\
		\forall j_1\in\Inink{l_1},\
		\forall p_1\in\Zn,\
		\forall l_2\in\Insys,\
		\forall j_2\in\Inink{l_2},\
		\forall p_2\in\Zn,\
		\hspace{1.7cm}
		\phantom{a}
	} \padDown
	\\
	\left(
		\bb{\mathcal{C}}^{[N]}
	\right)_{
		(j_1 + \sum_{\lambda=1}^{l_1-1} \nink{\lambda}),
		p_1,
		(j_2 + \sum_{\lambda=1}^{l_2-1} \nink{\lambda}),
		p_2
	} =
	\left\{
		\begin{array}{ll}
			\left(\mathcal{C}_{l_1}^{[N]}\right)_{j_1, p_1, j_2, p_2} & \hspace{-2mm} \text{ if } (l_1, j_1) = (l_2, j_2) \\
			0 & \hspace{-2mm} \text{ otherwise}
		\end{array}
	\right.
	\\
\end{array}
\end{equation}

Finally, \begin{small}$\matDUbar{1.5}{\bb{\hat{y}_C}^{[N+1]}}{\bb{\hat{\dot{y}}_C}^{[N+1]}}$\end{small} is a column vector of size $2 \nouttot$. It corresponds to a concatenation of the control parts in the COSTARICA estimators introduced in \ref{subsection:COSTARICA_estimator}.

\vspace{-8mm} %layout

\begin{equation}
\label{eq:GlobOp_yC_and_dyC}
\renewcommand{\arraystretch}{1.9}
\begin{array}{lccccr}
	\bb{\hat{y}_C}^{[N+1]} =
	\Big(
		&
		\hyC_{1, 1}^{[N+1]}\ ,\ ...\ ,\ \hyC_{1, \noutk{1}}^{[N+1]}\ ,
		\hsd & \hsd
		\hyC_{2, 1}^{[N+1]}\ ,\ ...\ ,\ \hyC_{2, \noutk{2}}^{[N+1]}\ ,
		\hsd & \hsd
		...\ ,
		\hsd & \hsd
		\hyC_{\nsys, 1}^{[N+1]}\ ,\ ...\ ,\ \hyC_{\nsys, \noutk{\nsys}}^{[N+1]}
		&
	\Big)^T
	\\
	\bb{\hat{\dot{y}}_C}^{[N+1]} =
	\Big(
		&
		\hdyC_{1, 1}^{[N+1]}\ ,\ ...\ ,\ \hdyC_{1, \noutk{1}}^{[N+1]}\ ,
		\hsd & \hsd
		\hdyC_{2, 1}^{[N+1]}\ ,\ ...\ ,\ \hdyC_{2, \noutk{2}}^{[N+1]}\ ,
		\hsd & \hsd
		...\ ,
		\hsd & \hsd
		\hdyC_{\nsys, 1}^{[N+1]}\ ,\ ...\ ,\ \hdyC_{\nsys, \noutk{\nsys}}^{[N+1]}
		&
	\Big)^T
	\\
\end{array}
\renewcommand{\arraystretch}{1.0}
\end{equation}

\vspace{-2mm} %layout

At this point, every quantity appearing in \eqref{eq:eta_MISSILES_N_solution_outputs_developed} has been defined. By manipulating the problem \eqref{eq:eta_MISSILES_N_solution_outputs_developed}, we can finally obtain the linear problem \eqref{eq:MISSILES}.

\vspace{-5mm} %layout

\begin{small}
\begin{equation}
\label{eq:MISSILES}
\begin{array}{||c||}
\hline\hline \\ \phantom{a}
	\left(
		\vphantom{\begin{array}{c}a\\a\\a\end{array}}
	 	I
	 	-
		\matDUbar{1.5}{\GbV^{[N]}}{\GbD^{[N]}}
		\bb{\mathcal{C}}^{[N]}
		\bb{\mathcal{A}}^{[N]}
		\matDDbar{1.5}{\!\!\!\Phi^T\!\!\!}{0}{0}{\!\!\!\Phi^T\!\!\!}
	\right)
	\matDUbar{1.5}{\hat{\yb}^{[N+1]}}{\hat{\dot{\yb}}^{[N+1]}}
	=
	\matDUbar{1.5}{\GbV^{[N]}}{\GbD^{[N]}}
	\bb{\mathcal{C}}^{[N]}
	\bb{\mathcal{B}}^{[N]}
	\!+\!
	\matDUbar{1.5}{\PbV^{[N]}}{\PbD^{[N]}}
	\bb{x}^{[N]}
	\!+\!
	\matDUbar{1.5}{\RbV^{[N]}}{\RbD^{[N]}}
	\bb{\tfC}^{[N]}
	\!+\!
	\matDUbar{1.5}{\bb{\hat{y}_C}^{[N+1]}}{\bb{\hat{\dot{y}}_C}^{[N+1]}}
\phantom{a} \\ \\ \hline\hline
\end{array}
\end{equation}
\end{small}

\vspace{-2mm} %layout

Problem \eqref{eq:MISSILES} can be resolved when all systems reached $t^{[N]}$, and gives a solution at $t^{[N+1]}$ by dispatching the solution of \eqref{eq:MISSILES} using \eqref{eq:MISSILES_solution_dispatching}. This solves the "Find root of this application" stage in figure \ref{fig:MISSILES_method}, and answers the the final question of subsection \ref{subsection:general_idea}.

\subsection{Implementation and first step}
\label{subsection:implementation_and_first_step}

In practice, it is possible to do the Hermite interpolation and the time-shift at the same time. Indeed, instead of computing $\bb{\mathcal{C}}^{[N]}$, $\bb{\mathcal{A}}^{[N]}$ and $\bb{\mathcal{B}}^{[N]}$ operators, it is possible to compute $\bb{\mathcal{C}}^{[N]}\bb{\mathcal{A}}^{[N]}$ and $\bb{\mathcal{C}}^{[N]}\bb{\mathcal{B}}^{[N]}$ directly to assemble the linear problem \eqref{eq:MISSILES_solution_dispatching}.

The underlying meaning of this replacement is that, instead of computing an interpolation on $t^{[N]}$ and $t^{[N+1]}$ and shifting it on $0$ and $\dt^{[N]}$, the interpolation is directly done on $0$ and $\dt^{[N]}$ with the same values and derivatives constraints. Indeed, no correction is required as $t^{[N+1]}-t^{[N]} = \dt^{[N]} - 0$.

Practically, a way to implement consists in:

\vspace{-0.5em} %layout

\begin{itemize}
\setlength{\itemsep}{-0.4em} %layout
\ite Replacing the $\bb{\mathcal{C}}^{[N]}$ operator by the identity tensor of order $4$, or equivalenty simply remove it from problem \eqref{eq:MISSILES}, and
\ite Computing $\bb{\mathcal{A}}^{[N]}$ and $\bb{\mathcal{B}}^{[N]}$ operators as explained in \eqref{eq:GlobOp_AVelem_ADelem_Belem_on_tN_tNp1} \eqref{eq:GlobOp_Acal} \eqref{eq:GlobOp_Bcal}, but replacing $t^{[N]}$ and $t^{[N+1]}$ by $0$ and $\dt^{[N]}$ respectively. Let's denote by $\AcVelem$, $\AcDelem$ and $\Bcelem$ the elementary quantities \eqref{eq:Hermite_AVelem_ADelem_Belem} with this change. Their (simpler) expressions in that case is given in \eqref{eq:AcVelem_AcDelem_Bcelem}.
\end{itemize}

\vspace{-6mm} %layout

\begin{equation}
\label{eq:AcVelem_AcDelem_Bcelem}
\renewcommand{\arraystretch}{2.2}
\begin{array}{c}
	\AcVelem =
	\left(0,\ 0,\ \dspfrac{3}{(\dt^{[N]})^2},\ \dspfrac{-2}{(\dt^{[N]})^3}\right)^T
	,\h
	\AcDelem =
	\left(0,\ 0,\ \dspfrac{-1}{\dt^{[N]}},\ \dspfrac{1}{(\dt^{[N]})^2}\right)^T
	,
	\\
	\Bcelem =
	\left(v_1,\ \dot{v}_1,\ \dspfrac{-3v_1}{(\dt^{[N]})^2}-2\dspfrac{\dot{v}_1}{\dt^{[N]}},\ \dspfrac{\dot{v}_1}{(\dt^{[N]})^2}+2\dspfrac{v_1}{(\dt^{[N]})^3}\right)^T
\end{array}
\renewcommand{\arraystretch}{1.0}
\end{equation}

\vspace{-1mm} %layout

Regarding the first macro-step, despite it is expected to know the initial values of all coupling variables (outputs and their corresponding inputs), there is usually no available time-derivatives of these coupling variables at $t^{[0]}=\tinit$. This makes it impossible to compute $\Bcelem$ and thus $\Belem^{[N]}$ in \eqref{eq:GlobOp_AVelem_ADelem_Belem_on_tN_tNp1} ($\dot{v}_1$ is not available).

In this case, the problem \eqref{eq:MISSILES} can still be assembled, but the underlying Hermite polynomials will simply be calibrated on three constraints: the input values at the beginning of the macro-step, and the input values and derivatives end the end of the macro-step. Analogously to the computations of subsection \ref{subsection:hermite_interpolation}, and calibrating the polynomials on $0$ and $\dt^{[N]}$ as explained above in this subsection, we obtain the expressions of $\AcVelemfs$, $\AcDelemfs$ and $\Bcelemfs$ in \eqref{eq:AcVelem_AcDelem_Bcelem_first_step}.

\vspace{-6mm} %layout

\begin{equation}
\label{eq:AcVelem_AcDelem_Bcelem_first_step}
\AcVelemfs \!=\!
\left(0,\ \frac{2}{\dt^{[N]}},\ \frac{-1}{(\dt^{[N]})^2}\right)^T
\!\!\!\!,
\hsd
\AcDelemfs \!=\!
\left(0,\ -1,\ \dspfrac{1}{\dt^{[N]}}\right)^T
\!\!\!\!,
\hsd
\Bcelemfs \!=\!
\left(v_1,\ \dspfrac{-2v_1}{\dt^{[N]}},\ \dspfrac{v_1}{(\dt^{[N]})^2}\right)^T
\end{equation}

\vspace{-2.5mm} %layout

\section{Results on test cases}
\label{section:results_on_test_cases}

This section presents results on benchmark co-simulation test cases. The MISSILES method is compared to the explicit fixed-step zero-order hold co-simulation method, also called non-iterative Jacobi (referred to as "\textbf{NI Jacobi}" in this section), as the latter is the most simple one, requiring none of the advanced capabilities mentioned in \ref{subsection:Capabilities} and therefore widely used in the industry. Comparisons are also done with the IFOSMONDI-JFM method \cite{Eguillon2021IfosmondiJFM}. As the latter uses an iterative Newton-like method (jacobian-free), the convergence criterion might affect the performance of a co-simulation. This criterion, further described in \cite{Eguillon2021IfosmondiJFM}, is based on a parameter called $\varepsilon$. The smaller this $\varepsilon$ is, the less tolerant the iterative method is on the validation of a solution.

The cases presented below here have been implemented in a way that enables all required capabilities, including the rollback, so that comparisons between the co-simulation methods can be made. However, in practice, as most of the modelling and simulation platforms provide non-rollback capable systems, neither IFOSMONDI-JFM nor classical IFOSMONDI methods can be used. Moreover, in order to get error measurements, we dispose of a monolithic simulation for each test case. That is to say, the simulation referred to as \textbf{monolithic reference} denotes the simulation of the global model  on a single solver, without coupling. Such simulations will be used as reference in this section. Please note that such monolithic simulation cannot be done in practice as the need for co-simulation usually arises when black-boxed systems (that might come from various simulation and modelling platforms) are connected to one another.

Despite the MISSILES method can handle variable-step co-simulation (different values of $\dt^{[N]}$ across $N$ can be taken), the time-stepping strategy was not discussed in this paper. Hence, results will be compared on co-simulation with a fixed macro-step size. The size of the macro-steps will be denoted by $\dt_{\text{fixed}}$, and we will simply have:

\vspace{-3mm} %layout

\begin{equation}
\label{eq:fixed_step}
\forall N\in[\![0, \Nmax[\![,\ \dt^{[N]} = \dt_{\text{fixed}}
\end{equation}

\vspace{-2mm} %layout

\subsection{Linear mechanical benchmark}
\label{subsection:linear_mechanical_benchmark}

This model is made of $\nsys=2$ systems. A sketch of it is presented on figure \ref{fig:Test_case_two_masses_sketch}. It is a very common benchmark for co-simulation methods \cite{Busch2016} \cite{Eguillon2019Ifosmondi} \cite{Meyer2021}.

\vspace{-7mm} %layout

\begin{center}
\includegraphics[scale=0.3]{\figuresdir/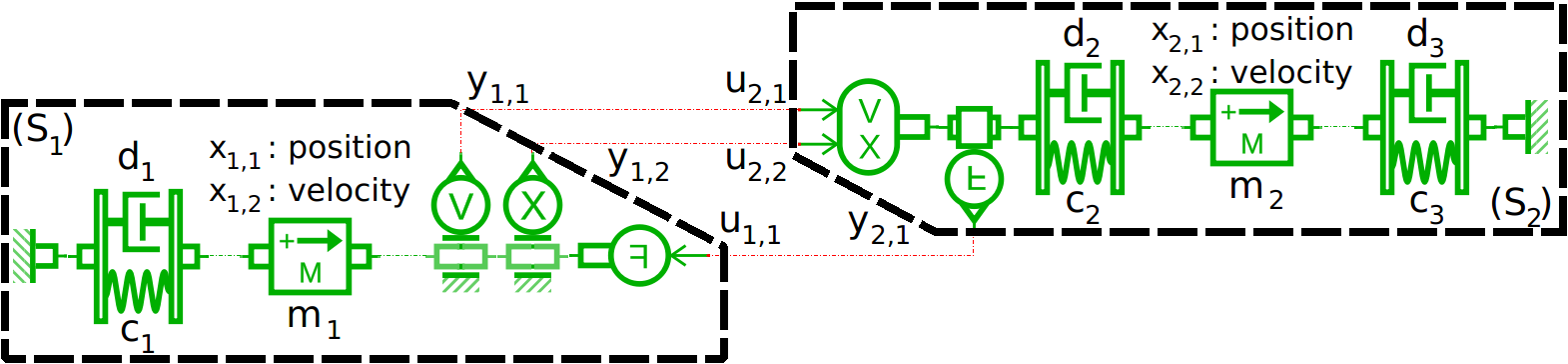}
\vspace{-3mm} %layout
\captionof{figure}{Benchmark co-simulation modular model: two linear mechanical bodies with springs and dampers}
\label{fig:Test_case_two_masses_sketch}
\end{center}

\vspace{-2mm} %layout

For the sake of reproducibility, the parameters of the bodies, springs, dampers and co-simulation run are given in table \ref{tab:Test_case_two_masses_parameters}. Physical quantities are measured positively from left to right, and negatively from right to left.

\vspace{-5mm} %layout

\begin{center}
\captionof{table}{Parameters of the linear mechanical test case}
\vspace{-1mm} %layout
\label{tab:Test_case_two_masses_parameters}
\begin{tabular}{cc}
	\renewcommand{\arraystretch}{1.0}
	\begin{tabular}{|c|c|c|}
		\hline
		\multicolumn{3}{|c|}{\textbf{Physical parameters}} \\
		\hline\hline
		\underline{Name} & \underline{Definition} & \underline{Value} \\
		$d_1$ & Damper rating & $10$ N/(m/s) \\
		$d_2$ & Damper rating & $10$ N/(m/s) \\
		$d_3$ & Damper rating & $40$ N/(m/s) \\
		$c_1$ & Sprint rate & $10\ 000$ N/m \\
		$c_2$ & Sprint rate & $10\ 000$ N/m \\
		$c_3$ & Sprint rate & $100\ 000$ N/m \\
		$m_1$ & Body mass & $5$ kg \\
		$m_2$ & Body mass & $80$ kg \\
		\hline
	\end{tabular}
	\renewcommand{\arraystretch}{1.0}
	&
	\begin{tabular}{r}
		\renewcommand{\arraystretch}{1.0}
		\begin{tabular}{|c|c|c|}
			\hline
			\multicolumn{3}{|c|}{\textbf{Initial states}} \\
			\hline\hline
			\underline{Expression} & \underline{Definition} & \underline{Value} \\
			$x_{1, 1}(\tinit)$ & Left body position & $-1$ m \\
			$x_{1, 2}(\tinit)$ & Left body velocity & $0$ m/s \\
			$x_{2, 1}(\tinit)$ & Right body position & $-3$ m \\
			$x_{2, 2}(\tinit)$ & Right body velocity & $0$ m/s \\
			\hline
		\end{tabular}
		\renewcommand{\arraystretch}{1.0}
		\padDown %layout
		\\
		\renewcommand{\arraystretch}{1.0}
		\begin{tabular}{|c|c|}
			\hline
			\multicolumn{2}{|c|}{\textbf{Co-simulation parameters}} \\
			\hline\hline
			\underline{Expression} & \underline{Value} \\
			$[\tinit, \tend]$ & $[0, 2]$ s \\
			\hline
		\end{tabular}
		\renewcommand{\arraystretch}{1.0}
	\end{tabular}
	\padDown %layout
	\\
	\multicolumn{2}{c}
	{
		\renewcommand{\arraystretch}{1.0}
		\begin{tabular}{|c|c|c|c|}
			\hline
			\multicolumn{4}{|c|}{\textbf{Initial coupling conditions}} \\
			\hline\hline
			\underline{Input} & \underline{Output} & \underline{Definition} & \underline{Value} \\
			$u_{1, 1}(\tinit)$ & $y_{2, 1}(\tinit)$ & Force on right of left mass & $-20\ 000$ N \\
			$u_{2, 1}(\tinit)$ & $y_{1, 1}(\tinit)$ & Left body velocity & $0$ m/s \\
			$u_{2, 2}(\tinit)$ & $y_{1, 2}(\tinit)$ & Left body position & $-1$ m \\
			\hline
		\end{tabular}
		\renewcommand{\arraystretch}{1.0}
	}
\end{tabular}
\end{center}

The results are presented in table \ref{tab:Test_case_two_masses_results}. The $\varepsilon$ parameter denotes the convergence criterion parameter of the iterative method used by IFOSMONDI-JFM, as described in the introduction of this section.

\vspace{-5mm} %layout

\begin{center}
\captionof{table}{Results on linear mechanical model: relative error on left body's position (in \%) and computational time (in s)}
\label{tab:Test_case_two_masses_results}
\renewcommand{\arraystretch}{1.5}
\begin{tabular}{|c||c|c|c|c|}
\hline
& \multirow{2}{*}{NI Jacobi} & \multicolumn{2}{|c|}{IFOSMONDI-JFM} & \multirow{2}{*}{MISSILES} \\
\cline{3-4}
& & $\varepsilon=10^{-2}$ & $\varepsilon=10^{-5}$ & \\
\hline\hline
$\dt_{\text{fixed}} = 10^{-3}$
	& \wordwrap{1.0}{$5.80$ \% \\ $0.26$ s}
	& \wordwrap{1.0}{$7.55\cdot 10^{-4}$ \% \\ $0.62$ s}
	& \wordwrap{1.0}{$2.66\cdot 10^{-4}$ \% \\ $0.95$ s}
	& \wordwrap{1.0}{$7.82\cdot 10^{-3}$ \% \\ $0.31$ s}
	\\
\hline
$\dt_{\text{fixed}} = 10^{-4}$
	& \wordwrap{1.0}{$2.93\cdot 10^{-1}$ \% \\ $1.80$ s}
	& \wordwrap{1.0}{$4.80\cdot 10^{-5}$ \% \\ $5.35$ s}
	& \wordwrap{1.0}{$2.27\cdot 10^{-5}$ \% \\ $7.48$ s}
	& \wordwrap{1.0}{$1.34\cdot 10^{-3}$ \% \\ $3.01$ s}
	\\
\hline
\end{tabular}
\renewcommand{\arraystretch}{1.0}
\end{center}

Several elements can be noticed on results of table \ref{tab:Test_case_two_masses_results}. The IFOSMONDI-JFM is slower than NI Jacobi and MISSILES as its iterative aspect make it require a larger amount of systems internal solvers restarts, which might be costly in terms of computational time.

Contrary to the IFOSMONDI-JFM method, when the accuracy is not satisfactory with MISSILES, there is no $\varepsilon$ parameter to tune to get a lower error for a given macro-step size, as the method is based on a direct solving of the coupling problem on each step. However, the macro-step size can be decreased in order to enhance the accuracy.

MISSILES is slower than NI Jacobi method (for a given fixed macro-step size) as the latter requires almost no computation in addition to the systems integrations. Please note that, in this case, as the macro-step size does not change and the systems $(S_1)$ and $(S_2)$ of figure \ref{fig:Test_case_two_masses_sketch} are linear, the matrix of the linear problem \eqref{eq:MISSILES} could be computed, assembled and factorized only once. However, in order to be as close as possible to a real use of co-simulation in practice, with black-boxed systems, we recomputed this operator at each macro-step, in order to mimic the case where the content of the systems is not known (modular models in industrial applications, for example). Nonetheless, this overhead in terms of computational time with respect to the NI Jacobi method is balanced by a better accuracy on MISSILES, as expected. 

Finally, despite the linear nature of the systems (making the COSTARICA estimators theoretically exact), the MISSILES method is not exact. Several causes can be mentioned: the inverse Laplace is done numerically with the Stehfest method \cite{Stehfest1970} \cite{Jacquot1983}, the global linear problem \eqref{eq:MISSILES} is solved numerically too (the accuracy is driven by the condition number of the matrix of this problem), the successive local polynomial approximations of non-polynomial solutions (coupling variables), ... Generally, MISSILES also relies on the data provided by the systems. For instance, the matrices of the linearizations are supposed to be exact, as their computation is done inside of the black-boxed systems. Nevertheless, the error reached by MISSILES on table \ref{tab:Test_case_two_masses_results} is satisfactory.

Indeed, figure \ref{fig:Test_case_two_masses_results} presents a superimposed view of the position of the left body ($x_{1, 1}$ state, also $y_{1, 2}$ output variable of system $(S_1)$) across the time, and it is noticable that the co-simulation with MISSILES is close to the monolithic reference quite as much as the co-simulation with the IFOSMONDI-JFM method. This is even more apparent on a zoom on a peak of this variable, as shown in figure \ref{fig:Test_case_two_masses_results_zoom}.

\vspace{-7mm} %layout

\begin{center}
\includegraphics[scale=0.5]{\figuresdir/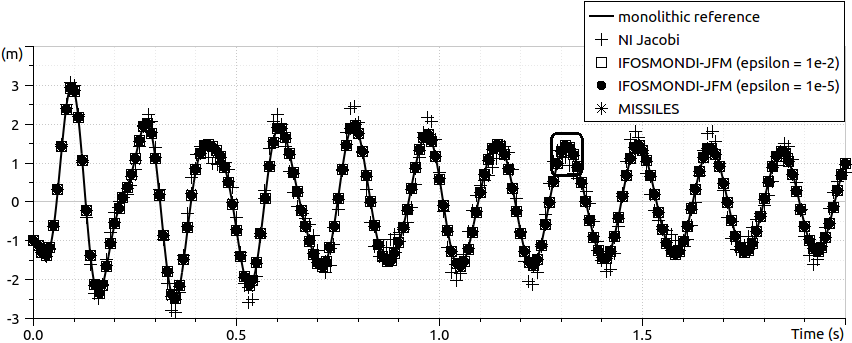}
\vspace{-3mm} %layout
\captionof{figure}{Left body's position - Comparison of co-simulation methods with $\dt_{\text{fixed}} = 10^{-3}$ s - The framed zone is zoomed on figure \ref{fig:Test_case_two_masses_results_zoom}}
\label{fig:Test_case_two_masses_results}
\end{center}

\vspace{-3mm} %layout

\begin{center}
\includegraphics[scale=0.5]{\figuresdir/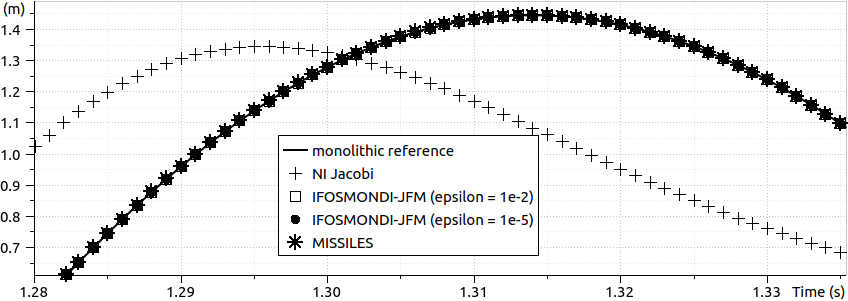}
\vspace{-3mm} %layout
\captionof{figure}{Zoom on $[1.29, 1.33]$ on curves of figure \ref{fig:Test_case_two_masses_results}}
\label{fig:Test_case_two_masses_results_zoom}
\end{center}

\vspace{-3mm} %layout

This proves the usefullness of the MISSILES method on such model, keeping in mind that the IFOSMONDI-JFM cannot be used in case all involved systems are not rollback-capable. MISSILES is a good way to reach an almost-similar accuracy, as shows figure \ref{fig:Test_case_two_masses_results_zoom}.

\subsection{Non-linear model: Lotka-Volterra equations}
\label{subsection:non_linear_model_lotka_volterra_equations}

This model is also made of $\nsys=2$ systems. A sketch of it is presented on figure \ref{fig:Test_case_LV_sketch}. The Lotka-Volterra prey-predator equations \cite{Volterra1928} are represented in this model, each system representing a species.

\vspace{-2mm} %layout

\begin{center}
\includegraphics[scale=0.3]{\figuresdir/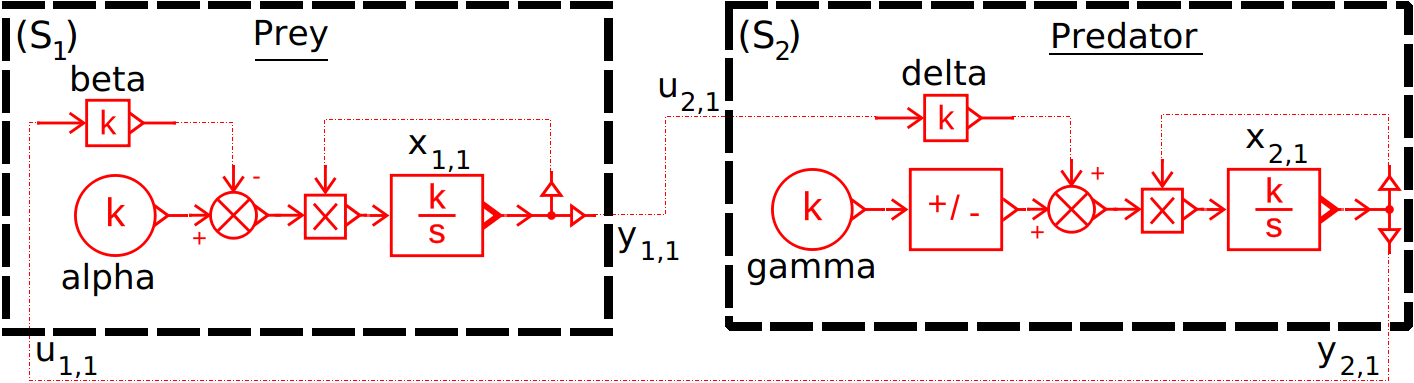}
\vspace{-3mm} %layout
\captionof{figure}{Lotka-Volterra model with a distinct model for the prey and another for the predator, in a co-simulation configuration}
\label{fig:Test_case_LV_sketch}
\end{center}

As for test case \ref{subsection:linear_mechanical_benchmark}, the parameters of the model are given in table \ref{tab:Test_case_LV_parameters} for the sake of reproducibility. The results are then presented in table \ref{tab:Test_case_LV_results}.

\vspace{-6mm} %layout

\begin{center}
\captionof{table}{Parameters of the Lotka-Volterra test case}
\vspace{-1mm} %layout
\label{tab:Test_case_LV_parameters}
\begin{tabular}{cc}
	\renewcommand{\arraystretch}{1.0}
	\begin{tabular}{|c|c|c|}
		\hline
		\multicolumn{3}{|c|}{\textbf{Model's parameters}} \\
		\hline\hline
		\underline{Name} & \underline{Definition} & \underline{Value} \\
		$\alpha$ & Natural prey's birth rate & $\nicefrac{2}{3}$ \\
		$\beta$ & Rate of predation upon the prey & $\nicefrac{4}{3}$ \\
		$\gamma$ & Predator-s natural death rate & $1$ \\
		$\delta$ & \wordwrap{1.0}{Growth rate upon predators \\(due to predation)} & $1$ \\
		\hline
	\end{tabular}
	\renewcommand{\arraystretch}{1.0}
	&
	\begin{tabular}{r}
		\renewcommand{\arraystretch}{1.0}
		\begin{tabular}{|c|c|c|}
			\hline
			\multicolumn{3}{|c|}{\textbf{Initial states}} \\
			\hline\hline
			\underline{Expression} & \underline{Definition} & \underline{Value} \\
			$x_{1, 1}(\tinit)$ & Amount of prey & $1$ \\
			$x_{2, 1}(\tinit)$ & Amount of predator & $1$ \\
			\hline
		\end{tabular}
		\renewcommand{\arraystretch}{1.0}
		\padDown %layout
		\\
		\renewcommand{\arraystretch}{1.0}
		\begin{tabular}{|c|c|}
			\hline
			\multicolumn{2}{|c|}{\textbf{Co-simulation parameters}} \\
			\hline\hline
			\underline{Expression} & \underline{Value} \\
			$[\tinit, \tend]$ & $[0, 20]$ s \\
			\hline
		\end{tabular}
		\renewcommand{\arraystretch}{1.0}
	\end{tabular}
	\padDown %layout
	\\
	\multicolumn{2}{c}
	{
		\renewcommand{\arraystretch}{1.0}
		\begin{tabular}{|c|c|c|c|}
			\hline
			\multicolumn{4}{|c|}{\textbf{Initial coupling conditions}} \\
			\hline\hline
			\underline{Input} & \underline{Output} & \underline{Definition} & \underline{Value} \\
			$u_{1, 1}(\tinit)$ & $y_{2, 1}(\tinit)$ & Amount of predator & $1$ \\
			$u_{2, 1}(\tinit)$ & $y_{1, 1}(\tinit)$ & Amount of prey & $1$ \\
			\hline
		\end{tabular}
		\renewcommand{\arraystretch}{1.0}
	}
\end{tabular}
\end{center}

As well as in the first test case, we observe a smaller error with the MISSILES method than with the NI Jacobi method for similar macro-step sizes. IFOSMONDI-JFM stays the method with the highest accuracy with $\dt_{\text{fixed}} = 10^{-3}$, which is not surprising as the model is non-linear. Indeed, the rollback is the only way to exactely solve the non-linear problem \eqref{eq:IFOSMONDIJFM_coupling_constraint} with \eqref{eq:IFOSMONDIJFM_etaN}. The COSTARICA estimators used by MISSILES are only locals and can thus only approximate the behavior of the systems on each step. This is namely the reason why, for smaller steps ($\dt_{\text{fixed}} = 10^{-4}$), the MISSILES method reaches a similar accuracy than the IFOSMONDI-JFM method (using real rollback). These results are presented in table \ref{tab:Test_case_LV_results}.

\vspace{-5mm} %layout

\begin{center}
\captionof{table}{Results on Lotka-Volterra model: relative error on prey (in \%) and computational time (in s)}
\label{tab:Test_case_LV_results}
\renewcommand{\arraystretch}{1.5}
\begin{tabular}{|c||c|c|c|c|}
\hline
& \multirow{2}{*}{NI Jacobi} & \multicolumn{2}{|c|}{IFOSMONDI-JFM} & \multirow{2}{*}{MISSILES} \\
\cline{3-4}
& & $\varepsilon=10^{-2}$ & $\varepsilon=10^{-5}$ & \\
\hline\hline
$\dt_{\text{fixed}} = 10^{-3}$
	& \wordwrap{1.0}{$1.34\cdot 10^{-1}$ \% \\ $1.01$ s}
	& \wordwrap{1.0}{$6.86\cdot 10^{-4}$ \% \\ $2.63$ s}
	& \wordwrap{1.0}{$6.86\cdot 10^{-4}$ \% \\ $3.39$ s}
	& \wordwrap{1.0}{$6.79\cdot 10^{-3}$ \% \\ $2.20$ s}
	\\
\hline
$\dt_{\text{fixed}} = 10^{-4}$
	& \wordwrap{1.0}{$1.35\cdot 10^{-2}$ \% \\ $13.27$ s}
	& \wordwrap{1.0}{$6.82\cdot 10^{-4}$ \% \\ $21.43$ s}
	& \wordwrap{1.0}{$6.82\cdot 10^{-4}$ \% \\ $22.60$ s}
	& \wordwrap{1.0}{$3.48\cdot 10^{-4}$ \% \\ $18.52$ s}
	\\
\hline
\end{tabular}
\renewcommand{\arraystretch}{1.0}
\end{center}

Regarding the comparison between the NI-Jacobi and the MISSILES method, an overshoot phenomenon cannot be seen on table \ref{tab:Test_case_LV_results} but can be observed on the solutions. Figure \ref{fig:Test_case_LV_results} shows the amount of prey ($x_{1, 1}$ state, also $y_{1, 1}$ output variable of system $(S_1)$). Macroscopically, all curves are superimposed, yet on the zoom on a peak presented in figure \ref{fig:Test_case_LV_results_zoom} we can observe that the NI Jacobi co-simulation produces an overshoot on the solution, where the MISSILES co-simulation produces results that are significantly more accurate and close to the monolithic reference solution.

The overshoot is a dangerous behavior in practice as it adds energy into the system. This is even more obvious on even greater macro-steps, if we look as the orbit of the solution. Such orbits are presented in figure \ref{fig:Test_case_LV_orbit_dt_1em2}, comparing the orbit of the solution on the monolithic simulation with the orbits with the NI Jacobi and the MISSILES method with a macro-step size of $0.01$ s. It can be notices on the plot in the middle that the NI Jacobi produced a solution with an diverging orbit, which is not the case with the MISSILES method.

\begin{center}
\includegraphics[scale=0.5]{\figuresdir/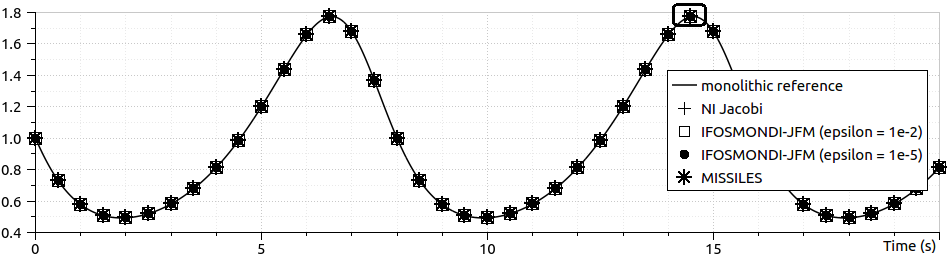}
\vspace{-3mm} %layout
\captionof{figure}{Amount of prey - Comparison of co-simulation methods with $\dt_{\text{fixed}} = 10^{-3}$ s - The framed zone is zoomed on figure \ref{fig:Test_case_LV_results_zoom}}
\label{fig:Test_case_LV_results}
\end{center}

\begin{center}
\includegraphics[scale=0.48]{\figuresdir/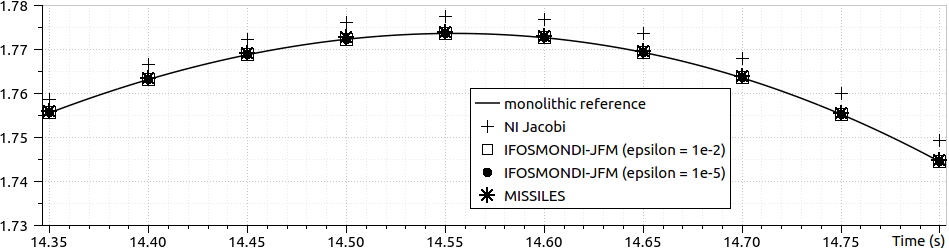}
\vspace{-3mm} %layout
\captionof{figure}{Zoom on $[14.35, 14.80]$ on curves of figure \ref{fig:Test_case_LV_results}}
\label{fig:Test_case_LV_results_zoom}
\end{center}

\vspace{-3mm} %layout

\begin{center}
\includegraphics[scale=0.58]{\figuresdir/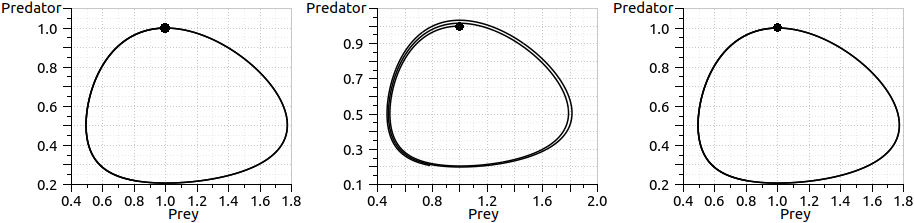}
\vspace{-4mm} %layout
\captionof{figure}{Orbit of solution on (from left to right): the monolithic reference, a co-simulation with $\dt_{\text{fixed}}=10^{-2}$ s with NI Jacobi method, a co-simulation with $\dt_{\text{fixed}}=10^{-2}$ s too and with MISSILES method - The initial state $(1, 1)$ is indicated by a black dot}
\label{fig:Test_case_LV_orbit_dt_1em2}
\end{center}

\vspace{-4mm} %layout

\section{Conclusion}
\label{section:conclusion}

The introduced MISSILES co-simulation method manages to properly approximate the implicit coupling constraint on the coupling variables between systems that are not capable of rollback, although the latter is a mandatory capability for the implicit coupling formulations' resolution methods. This approach makes it possible to satisfy an approximation of the whole set of contraints on all systems through the resolution of a single global system of linear equations. The coefficients of the polynomial expression of all coupling quantities can then be computed from the solution to this problem.

On linear co-simulation systems, this approach reaches a good time / accuracy trade-off: the computational time stays competitive with the classical zero-order hold non-iterative Jacobi fully explicit method due to the avoidance of repeated solver restarts (contrary to iterative co-simulation methods), an the accuracy stays competitive with implicit methods (requiring the rollback capability) as the solved constraint concerns the coupling variables at the end of the co-simulation steps.

Regarding non-linear systems, the MISSILES methods stays a satisfactory co-simulation method for large co-simulation steps in case not all systems are capable of rollback. Even if taking the non-linearities into account brings a significant improvement to the quality of the results in the context of an implicit co-simulation method (like it is the case in the rollback-based IFOSMONDI-JFM method), it is usually not possible to do it due to the scarsity of the rollback in practice. The accuracy of MISSILES on non-linear cases is related to the validity of the local linearizations of the systems. For this reason, small enough co-simulation steps can make the MISSILES method reach a accuracy that is competitive with the iterative IFOSMONDI-JFM co-simulation method. Depending on the case, the convenient co-simulation step size can be different, the aim being to stay on intervals where the linearization stays close enough to the systems' behaviors. Even in a given modular model, the satisfactory macro-step size might vary across the simulation time. For these reasons, the MISSILES method would benefit from an adaptive co-simulation time-stepper. Moreover, this would not be an obstacle to the construction of the core global system of linear equations of the method. Indeed, the whole paper was written without supposing a constant macro-step size.

In addition to that, as MISSILES provides an estimation of the coupling variables at the end of a macro-step before integrating the systems on it, this estimation can act as a predictor, a corrector being the coupling variables' values once the systems reach this time. Further research will investigate the usage of MISSILES' estimations in a time-stepper that could benefit the method.

\bibliographystyle{elsarticle-num}
\bibliography{MISSILES_refs}

\end{document}